\documentclass[12pt]{amsart}
\usepackage{amsfonts,amsmath,amssymb,comment,amscd,latexsym,{xy}}

\def\A{\mathbb{A}(\tl)}
\def\At{\mathbb{A}}
\def\c{\text{cc}}
\def\g{\mathcal{G}}
\def\cC{\mathcal{C}}
\def\E{\mathrm{Ext}_{[\tl]}(\k C_p,\k^G)}
\def\Z{\mathbb{Z}}
\def\X{\mathbb{X}(\tl)}
\def\x{\mathbb{X}}
\def\H{H^2_c(\k C_p,\k^G,\tl)}
\def\Hc{H^2_{\c}(\k C_p,\k^G,\tl)}
\def\t{\text{triv}}
\def\F{\Phi}

\def\we{\wedge}
\def\Tr{\text{tr}}

\newtheorem{Lem}{Lemma}[section]

\newtheorem{Prop}[Lem]{Proposition}
\newtheorem{Cor}[Lem]{Corollary}
\newtheorem{Thm}[Lem]{Theorem}

\theoremstyle{definition}
\newtheorem{Def}[Lem]{Definition}

\renewcommand{\k}{{\Bbbk}}
\renewcommand{\o}{\otimes}
\newcommand{\D}[1]{{\Bbbk}^{#1}}
\newcommand{\pf}{\medskip\noindent{\sc Proof:\,}}
\newcommand{\G}{\widehat G}
\newcommand{\w}[1]{\widehat #1}

\newcommand{\ov}[1]{\overline {#1}}
\newcommand{\bu}[1]{{#1}^{\bullet}}

\newcommand{\bi}{\begin{itemize}}
\newcommand{\ei}{\end{itemize}}

\newcommand\gen[1]{{\langle #1\rangle}}
\newcommand{\un}[1]{\underline #1}
\newcommand{\tl}{\triangleleft}
\newcommand{\tr}{\triangleright}
\numberwithin{equation}{section} 
\begin{document}

\title[Isomorphism types of Hopf algebras] {Isomorphism types of Hopf algebras in a class of abelian extensions. I.}
\author{Leonid Krop}
\thanks{Research  partially supported by a grant from the College of Liberal Arts and Sciences at DePaul University}
\address{DePaul University\\Chicago, IL 60614}\email{lkrop@depaul.edu}

\begin{abstract}There is no systematic general procedure by which isomorphism classes of Hopf algebras that are extensions of $\k F$ by ${\k}^G$ can be found. We develop the general procedure for classification of isomorphism classes of Hopf algebras which are extensions of the group algebra $\k C_p$ by ${\k}^G$ where $C_p$ is a cyclic group of prime order $p$ and ${\k}^G$ is the Hopf algebra dual of $\k G$, $G$ a finite abelian $p$-group and $\k$ is an algebraically closed field of characteristic $0$. We apply the method to calculate the number of isoclasses of commutative extensions and certain extensions of this kind of dimension $\le p^4$.
\end{abstract}
\date{4/14/14}
\maketitle



\section{Introduction}\label{intro}

 There is no systematic general procedure by which isomorphism classes of Hopf algebras that are extensions of $\k H$ by ${\k}^G$ can be found. The purpose of this article is to fill this gap in case $H=C_p$ and $G$ is a finite abelian $p$-group for a prime $p$, and $\k$ is an algebraically closed field  of characteristic zero.

Let us agree to write $\mathrm{Ext}(\k C_p,\k^{G})$ for the set of all equivalence classes of extensions of $\k C_p$ by $\k^G$. Elements of $\mathrm{Ext}(\k C_p,\k^{G})$ possess two special features. Every algebra $H$ there is equivalent as extension to smash product $\k^G\#\k C_p$ with respect to a certain action of $C_p$ on $\k^G$, and $\k^{C_p}$ is central in the dual Hopf algebra $H^*$. The action of $C_p$ on $\k^G$ induces an action `$\tl$' of $C_p$ on $G$, the corresponding $\mathbb{Z}C_p$-module is denoted by $(G,\tl)$. In consequence, $H$ is determined up to equivalence by a pair $(\tau,\tl)$ where $\tau:\k C_p\to\k^G\o\k^G$ is a $2$-cococycle deforming the tensor product coalgebra structure of $\k^G\o\k F$. Abelian extensions with undeformed multiplication were studied by M.Mastnak \cite{Ma}. We adopt a version of his notation $H^2_c(\k C_p,\k^G,\tl)$ for the group of Hopf $2$-cocycles.  

The first major result is a structure theorem for the group 

\noindent$H^2_c(\k C_p,\k^G,\tl)$. It states that if $G$ is any finite abelian $p$-group with $p>2$, or a finite elementary $2$-group then there is a $C_p$-isomorphism 
\begin{equation}\label{Hopf2cohomology}H^2_c(\k C_p,\k^G,\tl)\simeq H^2(C_p,\G,\tl)\times H^2_N(G,\bu{\k})\end{equation}
where $\G$ is the dual group of $G$, $H^2(C_p,\G,\tl)$ is the second cohomology group of $C_p$ over $\G$ with respect to the action `$\tl$', and $H^2_N(G,\bu{\k})$ is the kernel in the Schur multiplier of $G$ of the norm mapping. We point out that formula \eqref{Hopf2cohomology} can be seen as a generalization of the Baer's exact sequence for the cohomology group $H^2(G,C_p)$ of central extensions of $G$ by $C_p$ \cite [p.34]{BT}. For, setting $\tl=\text{triv}$, the trivial action, we show (Section \ref{commutative}) that   $H^2_c(\k C_p,\k^G,\text{triv})$ coincides with $H^2(G,C_p)$ while $H^2(C_p,\G,\text{triv})$ and $H^2_N(G,\k^{\bullet})$ can be identified with $\mathrm{Ext}^1_{\mathbb{Z}}(G,C_p)$ and $\text{Hom}(\wedge^2G,C_p)$, respectively. Hence (\ref{Hopf2cohomology}) turnes into the splitting Baer's exact sequence
$$1\to\mathrm{Ext}^1_{\mathbb{Z}}(G,C_p)\to H^2(G,C_p)\to \text{Hom}(\wedge^2G,C_p)\to 1$$

The principal result of the paper is the construction of a bijection between the orbits of a certain group $\g(\tl)$ in $\H$ and isoclasses of extensions.
$\g(\tl)$ is built from two groups $\A$ and $A_p$. There $\A$ is the group of all $C_p$-linear automorphisms of $(G,\tl)$ and $A_p=\mathrm{Aut}(C_p)$. $A_p$ does not act on $G$, but it does act on the set of homomorphisms $\tl:C_p\to\mathrm{Aut}(G)$ via 
\begin{equation*} a\tl^{\alpha} x=a\tl\alpha(x),\,\alpha\in A_p,a\in G,x\in C_p.\end{equation*}
Let $I(\tl,\tl')$ be the set of automorphisms of $G$ intertwining actions $\tl$ and $\tl'$. $I(\tl,\tl')$ is an $\A$-set, in fact a single orbit of $\A$ in $\mathrm{Aut}(G)$. For every $\alpha\in A_p$ fix some $\lambda_{\alpha}$ in $I(\tl,\tl^{\alpha})$. The group $\g(\tl)$ is the subgroup of $\mathrm{Aut}(G)$ generated by $\A$ and all $\lambda_{\alpha}$.  In fact, $\g(\tl)$ is a crossed product of $\A$ with a subgroup of $A_p$. It transpires that $H^2_c(\k C_p,\k^G,\tl)$ is a $\g(\tl)$-module. The orbits of $\g(\tl)$ in $\H$ determine isotypes of extension in the following way. Let us denote by $[\tl]$ the set of all actions $\tl'$ isomorphic to $\tl^{\alpha}$ for some $\alpha$. We designate $\mathrm{Ext}_{[\tl]}(\k C_p,\k^G)$ to the set of all equivalence classes of extensions whose $C_p$-action belongs to $[\tl]$. We recall that by a fundamental result of D. Stefan \cite{St} the number of isomorphism types in any of our classes is finite. In the case at hand, we show that there is a bijection between isotypes of noncocommutative Hopf algebras in $\mathrm{Ext}_{[\tl]}(\k C_p,\k^G)$ and the orbits of $\g(\tl)$ in $\H$ not contained in the subgroup $\Hc$ parametrizing cocommutative extensions.

Assuming $G$ elementary abelian and $p$ odd we extend the bijection to all isoclasses in $\mathrm{Ext}_{[\tl]}(\k C_p,\k^G)$. This is done by showing that for cocommutative extensions $\g(\tl)$-orbits in $\Hc$ coincide with $\A$-orbits there, and furthermore their isoclasses in are in $1-1$ correspondence with the orbits of $\A$ in $\Hc$. 

The last part of the paper is devoted to explicit calculations of the orbit set in several cases. For concrete calculations the smaller space $\mathbb{X}(\tl):=H^2(C_p,\G,\tl)\times H^2_N(G,\k^{\bullet})$ is most convenient. $\mathbb{X}(\tl)$ gets its $\g(\tl)$ action by transport of action via isomorphism \eqref{Hopf2cohomology}. The action is component-wise for every odd $p$. Assuming $G$ an elementary $p$-group and $p$ odd we show that there are $\lfloor\frac{3n+2}{2}\rfloor$ orbits for $\tl=\text{triv}$. We also describe orbit sets for all actions on elementary $p$- group $G$ of order $p^3$. Lower order cases, viz. $|G|=p,p^2$ are known with $|G|=p$ and $|G|=p^2$ due to \cite{Mas2} and \cite{Mas1}, respectively.

The paper is organized in six sections. In Section 1 we review the necessary facts of the theory of abelian extension. In section 2 we prove formula \eqref{Hopf2cohomology} for the groups $H^2_c(\k C_p,\k^G,\tl)$. Section 3 contains the isomorphism and bijection theorems. In Sections 4 and 5 we determine the orbit sets for commutative extensions with $G$ elementary $p$-group, and all extensions with $G$ elementary $p$-group of order $\le p^3$, respectively, and compute the number of isoclasses. 
\subsection{Notation and Convention}
In addition to notation introduced in the Introduction we will use the following.

 $A^{\bullet}$ the group of units of a commutative ring $A$.

 $\Gamma^n$ direct product of $n$ copies of group $\Gamma$.
 
 $\text{Fun}(\Gamma,A^{\bullet})$ the group of all functions from $\Gamma$ to $A^{\bullet}$ with pointwise multiplication. We will identify groups
 $\text{Fun}(G^n,(\k^{F^m})^{\bullet}),\\\text{Fun}(F^m,(\k^{G^n})^{\bullet})$ and $\text{Fun}(G^n\times F^m,\bu{\k})$ via $f(\un{a})(\un{x})=f(\un{x})(\un{a})=f(\un{a},\un{x})$ where $\un{a}\in G^n,\un{x}\in F^m$.
 
$Z^2(\Gamma,\bu{A},\bullet)),B^2(\Gamma,\bu{A},\bullet)$ and $H^2(\Gamma,\bu{A},\bullet)$ are the group of $2$-cocycles, $2$-coboundaries, and the second degree cohomology group of $\Gamma$ over $\bu{A}$ with respect to an action $\bullet$ of $\Gamma$ on $A$ by ring automorphisms.

$\delta_{\Gamma}$ the differential of the standard cochain complex for cohomology of the triple $(\Gamma,\bu{A},\bullet)$ \cite[IV.5]{Mac}.

$\Z_{p^n}$ cyclic group of order $p^n$ additively written.

By abuse of notation we will often use the same symbol for an element of $Z^2(\Gamma,\bu{A},\tl)$ and its image in $H^2(\Gamma,\bu{A},\tl)$.

Throughout the paper we treat the terms $\Gamma$-module, $\Gamma$-linear, etc as synonymous to $\Z\Gamma$-module, $\Z\Gamma$-linear, etc.

\section{Background Review}\label{background}
\subsection{Extensions of Hopf Algebras}
Let $\k$ be a ground field. In this paper we are concerned with finite-dimensional Hopf algebras over $\k$. For a Hopf algebra $H$ we use the standard notation $H^+=\text{Ker}\epsilon$. Let $\pi:H\to K$ be a morphism of Hopf algebras. We let $H^{\text{co}\pi}$ and ${}^{\text{co}\pi}H$ denote subalgebras of right/left coinvariants \cite[3.4]{Mo}. We adopt H.-J. Schneider's definition of a Hopf algebra extension \cite{S1}. In our context it is stated as follows.

\begin{Def}\label{weakextension} A Hopf algebra $C$ is an extension of a Hopf algebra $B$ by a Hopf algebra $A$ if there is a sequence of Hopf mappings 
\begin{equation} A\overset\iota\rightarrowtail C\overset\pi\twoheadrightarrow B\tag{E}.\end{equation}
with $\iota$ monomorphism, $\pi$ epimorphism, $\iota(A)$ normal in $C$ and $\text{Ker}\pi=\iota(A)^+C$.
\end{Def}
We add some comments to the definition. By \cite[Remark 1.2]{S2} or \cite[3.4.3]{Mo} we have the fundamental fact that $\iota(A)=C^{\text{co}\pi}$. It follows that our definition coincides with the definition of extension in \cite{A}. Conversely, the equality $\iota(A)=C^{\text{co}\pi}$ is equivalent to the equality $\iota(A)={}^{\text{co}\pi}C$ \cite[4.19]{BCM}, and both of them imply $\iota(A)$ is normal \cite[4.13]{BCM}. Even more is true. Either condition $\iota(A)=C^{\text{co}\pi}$ or $\text{Ker}\pi=\iota(A)^+C$ renders the sequence (E) an extension. For details see \cite[3.3.1]{A}. 

\subsection{Abelian Extensions}\label{abelian}

We assume in what follows the ground field $\k$ to be an algebraically closed field of characteristic $0$ and $C$ to be a finite-dimensional Hopf algebra. An extension (E) is called {\em abelian} if $A$ is commutative and $B$ is cocommutative. It is well-known \cite[Theorem 1]{LR} and \cite[2.3.1]{Mo} that in this case $A={\k}^G$ and $B=\k F$ for some finite groups $G$ and $F$. Below we consider only extensions of this kind and we use the notation
\begin{equation}\k^G\overset\iota\rightarrowtail H\overset\pi\twoheadrightarrow \k F.\tag{A}\end{equation}
To simplify notation we will refer to the Hopf algebra $H$ in a sequence (A) as an extension of $\k F$ by $\k^G$.
Essential to the theory of abelian extensions is a result in \cite{NVO}, or  general theorems \cite[2.4]{S2}, \cite[3.5]{MD}, asserting $H$ is a crossed product of $\k F$ over ${\k}^G$\footnotemark\footnotetext{ A short independent proof for abelian extension is given in the Appendix}. The theorem entails existence of a mapping called section ( see, e.g. \cite[3.1.13]{A}) 
\begin{equation}\label{section}\chi:\k F\to H\end{equation}
giving rise to the crossed product structure on $H$. This means $H={\k}^G\chi(F)$ with the multiplication
\begin{align}(f\chi(x))(f'\chi(y))&=f(\chi(x)f'\chi^{-1}(x))\chi(x)\chi(y)\label{multH}\\
&=f(\chi(x)f'\chi^{-1}(x))[\chi(x)\chi(y)\chi^{-1}(xy)]\chi(xy)\nonumber\end{align}
\noindent for $f,f'\in{\k}^G,\,x,y\in F$.
The mapping $x\o f\mapsto x.f:=\chi(x)f\chi^{-1}(x)$ defines a module algebra action of $F$ on ${\k}^G$ and the function $\sigma:F\times F\to{\k}^G,\sigma(x,y)=\chi(x)\chi(y)\chi^{-1}(xy)$ is a left, normalized $2$-cocycle for that action \cite[7.2.3]{Mo}. We recall that definition of action is independent of the choice of section, see e.g. \cite[7.3.5]{Mo}. 

We consider the dual of the above action of $\k F$ on $\k^G$. For any finite group $G$ we identify $(\k^G)^*$ with $\k G$ by treating $r\in\k G$ as the functional $f\mapsto f(r),f\in\k^G$. By general principles the transpose of a left module action of $\k F$ on $\k^G$ is a right module coalgebra action, denoted $\tl$ of $\k F$ on $(\k^G)^*=\k G$. Under this action $a\tl x$ is that element of $G$ for which 
\begin{equation}\label{dualFaction} (a\tl x)(f):=f(a\tl x)=(x.f)(a),\text{for all}\,f\in\k^G,a\in G,x\in F.\end{equation}
This definition makes sense as $\Delta_{\k G}$ is a $\k F$-linear map, hence there holds $\Delta_{\k G}(a\tl x)=a\tl x\o a\tl x$, whence $a\tl x$ is a grouplike, hence in $G$. We note that, in general, $\tl$ is a permutation action on $G$. Let $\{p_a|a\in G\}$ be a basis of $\k^G$ dual to the basis $\{a|a\in G\}$ of $\k G$. One can see easily that in the basis $\{p_a\}$ the two actions are related by the formula
\begin{equation}\label{.vstriangleleft}x.p_a=p_{a\triangleleft x^{-1}}\end{equation}

Theory of extensions has a fundamental duality expressed by the fact that for each sequence $(E)$ its companion sequence
\begin{equation}B^*\overset{\pi^*}\rightarrowtail C^*\overset{\iota^*}\twoheadrightarrow A^*\tag{$E^*$}\end{equation}
is also an extension, see \cite[4.1]{B} or \cite[3.3.1]{A}. Since for any finite group $F$, $(\k F)^*={\k}^F$ and $({\k}^F)^*=\k F$, every diagram (A) induces a diagram
\begin{equation}{\k}^F\rightarrowtail H^*\twoheadrightarrow \k G\tag{$A^*$}\end{equation}
A crossed product structure on $H^*$ is effected by a section
\begin{equation}\label{section*}\omega:\k G\to H^*\end{equation}
We choose to write $H^*=\omega(G){\k}^F$ with the multiplication
\begin{equation}\label{multH*}(\omega(a)\beta)(\omega(b)\beta')=\omega(ab)\tau(a,b)(\beta.b)\beta',\end{equation}
where for $a,b\in G,\beta,\beta'\in {\k}^F$
\begin{align}\beta.b&=\omega^{-1}(b)\beta\omega(b),\text{and}\\
\tau(a,b)&=\omega^{-1}(ab)\omega(a)\omega(b).\end{align}
We note that $\tau:G\times G\to {\k}^F$ is a right, normalized $2$-cocycle for the action $\beta\o b\mapsto\beta.b$. As above the right action of $G$ on $\D{F}$ induces a left action of $G$ on $F$ by permutations denoted by $a\triangleright x$, and the two actions are related by
\begin{equation}\label{.vstriangleright}p_x.a=p_{a^{-1}\triangleright x}\end{equation}
We fuse both actions into the definition of a product on $F\times G$ via
\begin{equation}\label{grpmult}(xa)(yb)=x(a\triangleright y)(a\triangleleft y)b\end{equation}
It was noted by M.Takeuchi \cite{T} that the composition \ref{grpmult} defines a group structure on $F\times G$ provided the actions $\tr,\tl$ satisfy the conditions
\begin{align}\label{matchedactions} ab\tl x&=(a\tl(b\tr x))(b\tl x)\\
               a\tr xy&=(a\tr x)((a\tl x)\tr y)\end{align}
We use the standard notation $F\bowtie G$ for the set $F\times G$ endowed with multiplication \eqref{grpmult}. We will also adopt the notation $\ov{x}$ for $\chi(x)$ and $\ov{a}$ for $\omega(a)$.

The above discussion enables us to associate a datum $\{\sigma,\tau,\tl,\tr\}$ to every Hopf algebra $H$ in an extension of type (A), and we write $H=H(\sigma,\tau,\tl,\tr)$ and $H^*=H^*(\sigma,\tau,\tl,\tr)$ for $H$ and its dual.

\subsection{Cocentral Extensions}\label{cocentral}
An extension (A) is called cocentral \cite{KMM} if $\D{F}$ is a central subalgebra of $H^*$. We record two properties of cocentral extensions needed below
\begin{Lem}\label{cocentralext} {\rm (1)} An extension {\rm(A)} is cocentral iff $\tr$ is trivial, or equivalently $G$ is a normal subgroup of $F\bowtie G$ in which case $F\bowtie G$ is a semidirect product $F\ltimes G$.

{\rm(2)} If {\rm(A)} is cocentral, then $\Delta_{\k^G}$ is $F$-linear.
\end{Lem}
\pf (1)It is well known \cite[(4.10)]{M} that $F\bowtie G$ is a group for any two actions $\tl,\tr$ arizing from an abelian extension. By \eqref{grpmult} $x^{-1}ax=x^{-1}(a\tr x)(a\tl x)$, hence $x^{-1}ax\in G$ for all $a\in G,x\in F$ iff $x^{-1}(a\tr x)=1$, that is $a\tr x=x$. On the other hand $\D{F}$ is cocentral iff $p_x.a=p_{a^{-1}\tr x}=p_x$ by \eqref{.vstriangleright}. The rest of part (1) is immediate from \eqref{matchedactions}.

(2) We must show the equality $\Delta_{\k^G}(x.p_a)=x.\Delta_{\k^G}(p_a)$. On the one hand we have
\begin{equation*}x.\Delta_{\k^G}(p_a)=\sum_{bc=a} x.p_b\o x.p_c=\sum_{bc=a}p_{b\tl x^{-1}}\o p_{c\tl x^{-1}}.\end{equation*}
In the second place 
\begin{equation*}\Delta_{\k^G}(x.p_a)=\Delta_{\k^G}(p_{a\tl x^{-1}})=\sum_{ef=a\tl x^{-1}}p_e\o p_f\end{equation*}
It remains to notice that the mappings $(b,c)\mapsto (b\tl x^{-1},c\tl x^{-1})$,\\$(e,f)\mapsto (e\tl x,f\tl x)$ give a bijective correspondence between the sets $\{(b,c)|bc=a\}\;\text{and}\,\{(e,f)|ef=a\tl x^{-1}\}$ as the action `$\tl$' is by group automorphisms.\qed 

Below we will write an extension datum $\{\sigma,\tau,\tl\}$ when $\tr$ is trivial. We will need explicit formulas for coalgebra structure mappings on $H$ and $H^*$ expressed in terms of their datum. These are the duals of the algebra structures \eqref{multH} and \eqref{multH*}, and follow from \cite[(4.5)]{M}.
\begin{Prop}\label{cocrossedproduct1}Let $H$ and $H^*$ be defined by a datum $\{\sigma,\tau,\tl\}$. The coalgebra structure of $H$ and $H^*$ is given by the mappings
\begin{align}\Delta_H(f\ov{x})&=\sum_{a,b\in G}\tau(x,a,b)f_1p_a\ov{x}\o f_2p_b\ov{x},\label{DeltaH}\\
\epsilon_H(f\ov{x})&=f(1_G).\nonumber\end{align}
\begin{align}\label{Delta*}\Delta_{H^*}(\ov{a}g)&=\sum_{x,y\in F}\sigma(x,y,a)\ov{a} p_xg_1\o(\ov{a\tl x}) p_yg_2,\\
\epsilon_{H^*}(\ov{a}g)&=g(1_F),\nonumber\end{align}\end{Prop}
\noindent where $f\in \k^G,g\in\k^F$.

For discussion of cohomology of abelian extensions we introduce the subgroup $\mathrm{Map}(F^n\times G^m,\bu{\k})$ of $\text{Fun}(F^n\times G^m,\bu{\k})$ of normalized $n+m$-dimensional cochains, i.e. functions satisfying $f(x_1,\ldots,x_n,a_1,\ldots,a_m)\\=1$ if at least one component of $(x_1,\ldots,x_n,a_1,\ldots,a_m)$ is the identity. Suppose $F$ acts on $G$ via $\tl$. We extend this action to $\mathrm{Map}(F^n\times G^m,\bu{\k})$ by the rule
\begin{equation}\label{Faction}y.f(x_1,...,x_n,a_1,...,a_m)=f(x_1,...,x_n,a_1\tl y,\ldots, a_m\tl y).\end{equation}
Identifying $\mathrm{Map}(F^n\times G^m,\bu{\k})$ with either $\mathrm{Map}(F^n,\bu{(\k^{G^m})})$ or\\ $\mathrm{Map}(G^m,\bu{(\k^{F^n})})$ (see \cite{Ma}) we denote by $\delta_F$,$\delta_G$, respectively, the standard differentials of group cohomology.
Finally we state conditions for equivalence \cite[3.4]{M} of two extensions in the form needed below. They are a particular case of\cite[5.2]{M}. 
\begin{Lem}\label{equivalence} Two extensions $H$ and $H'$ defined by data $\{\sigma,\tau,\tl\}$ and $\{\sigma',\tau',{\tl}'\}$ are equivalent if and only if $\tl={\tl}'$ and there exists $\zeta\in\text{Map}(F\times G,\bu{\k})$ satisfying
\begin{equation}\label{equivsigma}\sigma'=\sigma\delta_F\zeta^{-1}\,\text{and}\;\tau'=\tau\delta_G\zeta\end{equation}
If so, an isomorphism $\psi:H\to H'$ defined by $\psi(f\ov{x})=f\zeta(x)\ov{x},f\in\k^G,x\in F$ carries out the equivalence. 
\end{Lem}

\subsection{Cohomology Groups $H^2_c(\k F,\k^G,\tl)$}\label{cohomologygroups}
We describe in some detail cohomology theory of abelian extensions in the special case studied below. In the rest of the paper we consider cocentral extensions (A) satisfying the condition
\begin{equation}\label{cextensions}H^2(F,\bu{(\k^G)},\tl)=\{1\}\, \text{for every action} \tl.\end{equation}
We observe a simple
\begin{Lem}\label{reduction} Under the assumption \eqref{cextensions} every extension $H(\sigma,\tau,\tl)$ is equivalent to an extension $H(1,\tau',\tl)$\end{Lem}
\pf The condition \eqref{cextensions} means $\sigma=\delta_{F}\zeta,\zeta\in\mathrm{Map}(F\times G,\bu{\k})$. Then by Lemma \ref{equivalence} the mapping $\psi:H(\sigma,\tau,\tl)\to H(1,\tau',\tl),\,\psi(f\ov{x})=f\zeta(x)\ov{x},f\in \D{G},x\in F$ with $\tau'=\tau\delta_G\zeta$ is the required equivalence.\qed

Below we will write $H=H(\tau,\tl)$ for an extension with a datum $\{\sigma,\tau,\tl\}$ and $\sigma$ trivial. We let $\mathrm{Ext}(\k F,\k^G)$ stand for the set of equivalence classes of extensions of type(A), and we write $\mathrm{Ext}(\k F,\k^G,\tl)$ for the group of equivalence classes of extensions with fixed action $\tl$.

Cocentral extensions satisfying \eqref{cextensions} have been studied in \cite{Ma}. It is shown there \cite[4.4]{Ma} that the standard second cohomology group of abelian extensions \cite{H} coincides with the one in the next definition.
\begin{Def}\label{Hopfcocycles}  We let $Z^2_c(\k F,\D{G},\tl)$ denote the subgroup of all elements $\tau$ of $Z^2(G,\bu{(\k^F)},\text{id})$ satisfying $\delta_{F}(\tau)=\epsilon$. We let $B^2_c(\k F,\D{G},\tl)$ stand for the subgroup of $2$-cocycles $\delta_{G}\eta,\,\eta\in\mathrm{Map}(F\times G,\bu{\k})$ satisfying $\delta_{F}\eta=1$. We define the second degree Hopf cohomology group by
\begin{equation*}H^2_c(\k F,\D{G},\tl)=Z^2_c(\k F,\D{G},\tl)/B^2_c(\k F,\D{G},\tl).\end{equation*}\end{Def}
Explicitly both conditions $\delta_{F}\tau=\epsilon$ and $\delta_{F}\eta=\epsilon$ are expressed by:
\begin{align}\tau(xy)&=\tau(x)(x.\tau(y))\label{hopfcocycle}\\\eta(xy)&=\eta(x)(x.\eta(y))\label{hopf1cocycle}\end{align}
for all $x,y\in F$. 

We call elements of $Z^2_c(\k F,\k^G,\tl)$ and $B^2_c(\k F,\k^G,\tl)$ {\em Hopf} $2$-cocycles and $2$-coboundaries, respectively. We will use abbreviated symbols $Z^2_c(\k F,\k^G,\tl),B^2_c(\k F,\k^G,\tl)$, etc. for $Z^2_c(\tl),B^2_c(\tl)$, etc when the groups $G$ and $F$ are fixed. The restriction $B^2(G,(\k^F)^{\bullet})\cap Z^2_c(\tl)$ of the group of coboundaries to $Z^2_c(\tl)$ will be denoted by $B^2_{\c}(\tl)$. As $B^2_c(\tl)\subset B^2_{\c}(\tl)$ we can form the subgroup $H^2_{\c}=B^2_{\c}(\tl)/B^2_c(\tl)$ of $H^2_c(\tl)$. Explicitly, 
\begin{equation*}H^2_{\c}(\k F,\k^G,\tl):=\{\delta_{G}\eta B^2_c(\tl)|\eta\in\mathrm{Map}(F\times G,\bu{\k})\, \text{with}\,\delta_F(\delta_{G}\eta)=1\} .\end{equation*} 
We need to establish $F$-invariance of subgroups just defined.
\begin{Lem}\label{Finvariance} If $F$ is abelian, then groups $Z^2_c(\tl),B^2_{\c}(\tl)$, and $B^2_c(\tl)$ are $F$-invariant.\end{Lem}

\pf Suffices to show that for any $f\in\mathrm{Map}(F\times G^m,\bu{\k})$ and $g\in\mathrm{Map}(F^n\times G,\bu{\k})$ there holds
\begin{align} x.\delta_{F} f&=\delta_{F}(x.f)\label{delta1},\\
x.\delta_{G}\,g&=\delta_{G}(x.g)\label{delta2}
\end{align}
For, setting $m=1,2,f=\eta,\tau$ and $n=1,g=\eta$ in \eqref{delta1} and \eqref{delta2}, respectively we get our statement.

\eqref{delta2} is immediate from definitions in view of $G$ acting trivially on $F$. To show \eqref{delta1} we calculate
\begin{align*} &(x.\delta_{F}f)(y,z,a_1,\ldots,a_m)=(\delta_{F}f)(y,z,a_1\tl x,\ldots,a_m\tl x)=\\
&f(y,a_1\tl x,\ldots,a_m\tl x)(y.f)(z,a_1\tl x,\ldots,a_m\tl x)\\&f(yz,a_1\tl x,\ldots,a_m\tl x)^{-1}
=f(y,a_1\tl x,\ldots,a_m\tl x)\\&(xy.f)(z,a_1x,\ldots,a_m)f(yz,a_1\tl x,\ldots,a_m\tl x)^{-1}
\end{align*}
Switching around $x$ and $y$ in the middle term we get exactly $\delta_{F}(x.f)$.\qed

\section{Structure of $H^2_c(\k C_p,\k^G,\tl)$}
Unless stated otherwise $G$ is a $p$-group and $F=C_p$. The group $C_p\bowtie G$ is a $p$-group as well, hence nilpotent. As $G$ has index $p$ in $C_p\bowtie G$, $G$ is normal in $C_p\bowtie G$. By Lemma \ref{cocentralext}(1) the action $\tr$ is trivial. In addition, as $\bu{\k}$ is a divisible group the group $H^2(C_p,\bu{(\k^G)},\tl)$ vanishes by e.g. \cite[4.4]{Ma}. Thus the results of the preceeding section are applicable. We note a simple fact.
\begin{Lem} Let $\tau\in Z^2(G,\bu{(\D{F})})$. Then for every $x\in C_p\,\tau(x)$ is a $2$-cocycle for $G$ with coefficients in $\bu{\k}$ with the trivial action of $G$ on $\bu{\k}$.\end{Lem}
\pf The $2$-cocycle condition for the trivial action  is 
\begin{equation}\label{2cocycle} \tau(a,bc)\tau(b,c)=\tau(ab,c)\tau(a,b).\end{equation} 
Expanding both sides of the above equality in the basis $\{p_x\}$ and equating coefficients of $p_x$ proves the assertion.\qed 

Consider group $F$ acting on an abelian group $A$, written multiplicatively, by group automorphisms. Let $\mathbb{Z}F$ be the group algebra of $F$ over $\mathbb{Z}$. $\mathbb{Z}F$ acts on $A$ via
\begin{equation*} (\sum c_ix_i).a=\prod x_i.(a^{c_i}),\,c_i\in\mathbb{Z},\,x_i\in F.
\end{equation*}
For $F=C_p$ pick a generator $t$ of $C_p$ and set $\phi_i=1+t+\cdots +t^{i-1}, i=1,\ldots,p$. Choose $\tau\in Z^2(G,\bu{(\D{C_p})})$ and expand $\tau$ in 
terms of the standard basis $p_{t^i}$ for $\k^{C_p}$, $\tau=\sum\tau(t^i)p_{t^i}$ with $\tau(t^i)\in Z^2(G,\bu{\k})$. An easy induction on $i$ shows that condition \eqref{hopfcocycle} implies  
\begin{equation}\label{componentsoftau} \tau(t^i)=\phi_i.\tau(t),\;\text {for all}\;i=1,\ldots,p\end{equation}
For $i=p$ we have
\begin{equation}\label{cyclotomic} \phi_p.\tau(t)=1\end{equation}
in view of $t^p=1$ and $\tau(1)=1$.

Let $M$ be a $\mathbb{Z}C_p$-module. Following \cite{Mac} we define the mapping $N: M\to M\,\text{by}\, N(m)=\phi_p(t).m$. We denote by $M_N$ the kernel of $N$ in $M$. For $M=Z^2(G,\bu{\k}), B^2(G,\bu{\k})$ or $H^2(G,\bu{\k})$ we write $Z^2_N(G,\bu{\k})$ for $Z^2(G,\bu{\k})_N$ and similarly for the other groups. We abbreviate $Z^2_N(G,\bu{\k})$ to $Z^2_N(\tl)$ and likewise for $B^2_N(G,\bu{\k})$ and $H^2_N(G,\bu{\k})$.
\begin{Def}\label{admissible} We call a $2$-cocycle $s\in Z^2(G,\bu{\k})$ {\em admissible} if $s$ satisfies the condition
\begin{equation}\label{Adm}\phi_p.s=1\end{equation}\end{Def}
Thus by definition the set of all admissible cocycles is $Z^2_N(\tl)$. We note that $Z^2_N(\tl)$ is a subgroup of $Z^2(G,\k^{\bullet})$ as $\mathbb{Z}C_p$ acts by endomorphisms of $Z^2(G,\k^{\bullet})$. We want to compare abelian groups $Z^2_c(\tl)$ and $Z^2_N(\tl)$. This is done via the mapping 
\begin{equation*}\Theta:Z^2(G,(\k^{C_p})^{\bullet})\to Z^2(G,\k^{\bullet}),\,\Theta(\tau)=\tau(t).\end{equation*}
\begin{Lem}\label{hf=adm} The mapping $\Theta$ induces a $C_p$-isomorphism between $Z^2_c(\tl)$ and $Z^2_N(\tl)$.\end{Lem}
\pf We begin with an obvious equality $x.(\tau(y))=(x.\tau)(y)$. Taking $y=t$ we get $\Theta(x.\tau)=x.\Theta(\tau)$, that is $C_p$-linearity of $\Theta$. The relations \eqref{componentsoftau} show that $\Theta$ is monic. It remains to establish that $\Theta$ is epic.

Suppose $s$ is an admissible $2$-cocycle of $G$ in $\bu{\k}$. Define $\tau: G\times G\to\bu{(\D{C_p})}$ by setting $\tau(t^i)=\phi_i(t).s,\,1\le i\le p$. The proof will be complete if we demonstrate that $\tau$ satisfies \eqref{hopfcocycle}. 

For any $i,j\le p$ we have
\begin{equation*} \tau(t^i)(t^i.\tau(t^j))=(\phi_i(t).s)(t^i\phi_j(t).s)=(\phi_i(t)+t^i\phi_j(t)).s\end{equation*}
One sees easily that $\phi_i(t)+t^i\phi_j(t)=\displaystyle\sum_{k=0}^{i+j-1}t^k$. Hence if  $i+j< p$ we have $\phi_i(t)+t^i\phi_j(t)=\phi_{i+j}(t)$ and so $\tau(t^i)(t^i.\tau(t^j))=\tau(t^{i+j})$. If $i+j=p+m$ with $m\ge 0$, then $\displaystyle\sum_{k=0}^{p+m-1}t^k=\phi_p(t)+t^p(1+\cdots+t^{m-1})$ which implies $(\displaystyle\sum_{k=0}^{p+m-1}t^k).s=\phi_p(t).s\cdot t^p\phi_m(t).s=\phi_m(t).s=\tau(t^{i+j})$ by \eqref{Adm} and as $t^p=1$.\qed

The next step is to describe structure of $H^2_{\c}(\tl)$. We need some preliminaries. First, we write $x.f$ for the left action of $C_p$ on $\k^G$ dual to $\tl$ as in \eqref{dualFaction}. Since $\G$ is the group of grouplikes of $\k^G,\,\G$ is $C_p$-stable by Lemma \ref{cocentralext}(2). Further, we use
$\delta$ for the differential on the group of $1$-cochains of $G$ in $\bu{\k}$. We also note  $B^2_N(\tl)=B^2(G,\bu{\k})\cap Z^2_N(\tl)$. By \eqref{Adm} $\delta f\in B^2_N(\tl)$ iff $\phi_p(t).\delta f=1$ which, in view of $\delta$ being $C_p$-linear, is the same as $\delta(\phi_p(t).f)=1$. Since $(\delta f)(a,b)=f(a)f(b)f(ab)^{-1}$, $\text{Ker}\,\delta$ consists of characters of $G$, whence $\delta f\in B^2_N(\tl)$ iff $\phi_p(t).f$ is a character of $G$. Say $\chi=\phi_p(t).f\in \widehat{G}$. Then as $t\phi_p(t)=\phi_p(t)$, $\chi$ is a fixed point of the $C_p$-module $\widehat{G}$. Letting $\widehat{G}^{C_p}$ stand for the set of fixed points in $\G$ we have by \cite[IV.7.1]{Mac} an isomorphism 
$H^2(C_p,\widehat{G},\bullet)\simeq\widehat{G}^{C_p}/N(\widehat{G})$.  We connect $B^2_N(\tl)$ to $H^2(C_p,\widehat{G})$ via the homomorphism 
\begin{equation}\label{coboundarymap} \Phi:B^2_N(\tl)\to H^2(C_p,\widehat{G},\bullet),\,\delta f\mapsto (\phi_p.f)N(\widehat{G})\end{equation}         \begin{Lem}\label{structureofcoboundaries} The following properties holds
\bi\item[(i)] $\Theta(B^2_{\c}(\tl))=B^2_N(\tl)$,
\item[(ii)] $\Theta(B^2_c(\tl))=\ker\Phi$,
\item[(iii)] $B^2_N(\tl)/\ker\Phi\simeq H^2(C_p,\widehat{G},\bullet)$,
\item[(iv)] $H^2_{\c}(\tl)\simeq H^2(C_p,\widehat{G},\bullet)$.
\ei \end{Lem}                                                                                                                                                \pf First we show that $\Phi$ is well-defined. For, $\delta f=\delta g$ iff $fg^{-1}=\chi\in\widehat{G}$, hence
\begin{align*}&\Phi(\delta f)=(\phi_p.f)N(\widehat{G})=(\phi_p.g\chi)N(\widehat{G})\\
&=(\phi_p.g\cdot\phi_p.\chi)N(\widehat{G})=(\phi_p.g)N(\widehat{G})=\Phi(\delta g)\end{align*}

(i) Take some $\delta_G\eta\in B^2_{\c}(\tl)$. Evidently for every $x\in C_p$ (*) $(\delta_G\eta)(x)=\delta(\eta(x))$, hence $\Theta(\delta_G\eta)=\delta(\eta(t))$ is a coboundary, and $\phi_p.\delta(\eta(t))=1$ by \eqref{cyclotomic}, whence $\Theta(\delta_G\eta)\in B^2_N(\tl)$.
Conversely, pick $\delta f\in B^2_N(\tl)$ and define $\omega=\sum_{i=1}^p(\phi_i.\delta f)p_{t^i}$. The argument of Lemma \ref{hf=adm} shows $\omega$ lies in $Z^2_c(\tl)$. Set $\eta=\sum_{i=1}^p(\phi_i.f)p_{t^i}$. Using (*) again we derive 
\begin{equation*}\delta_G\eta=\sum_{i=1}^p(\phi_i.\delta f)p_{t^i}=\omega, \end{equation*}
hence $\delta_G\eta\in B^2_{\c}(\tl)$. Clearly $\Theta(\delta_G\eta)=\delta f$.

(ii) The argument of Lemma \ref{hf=adm} is applicable to $1$-cocycles satisfying \eqref{hopf1cocycle}. It shows that $\eta$ satisfies \eqref{hopf1cocycle} iff
\begin{equation}\label{hopf1cocycle'}\eta(t^i)=\phi_i.\eta(t)\end{equation}                                                                                  For $i=p$ we get $\phi_p.\eta(t)=\epsilon$, hence the calculation
\begin{align*}\Phi(\Theta(\delta_G\eta))=\Phi(\delta(\eta(t)))=(\phi_p.\eta(t)) N(\widehat{G})=N(\widehat{G}).\end{align*}
gives one direction. Conversely, $\Phi(\delta f)\in N(\widehat{G})$ means $\phi_p.f=\phi_p.\chi$ which implies $\phi_p.f\chi^{-1}=\epsilon$. Set $g=f\chi^{-1}$ and define   $1$-cocycle $\eta_g=\sum_{i=1}^p(\phi_i.g)p_{t^i}$. Since $\phi_p.g=\epsilon$,  $\eta_g$ satisfies \eqref{hopf1cocycle}, whence  $\delta_G\eta_g\in B^2_c(\tl)$. As $(\delta_G\eta_g)(t)=\delta g=\delta f$ by construction, $\Theta(\delta_G\eta_g)=\delta f$.                                                                                   
            
(iii) We must show that $\Phi$ is onto. For every character $\chi$ in $\widehat{G}^{C_p}$ we want to construct an $f:G\to\bu{\k}$ satisfying $\phi_p.f=\chi$. To this end we consider splitting of $G$ into the orbits under the action of $C_p$. Since every orbit is either regular, or a fixed point we have
\begin{equation*} G=\cup_{i=1}^r\{g_i,g_i\tl t,\ldots,g_i\tl t^{p-1}\}\cup G^{C_p}\end{equation*} 
For every $s\in G^{C_p}$ we pick a $\rho_s\in\k$ satisfying $\rho^p_s=\chi(s)$. We define $f$ by the rule
\begin{align*}f(g_i)=\chi(g_i),\,f(g_i\tl t^j)&=1\;\text{for all}\;j\ne 1\,\text{and all}\;i=1,\ldots,r,\,\text{and}\\
f(s)&=\rho_s\;\text{for every}\;s\in G^{C_p}\end{align*}                                                                                                     By definition $(\phi_p.f)(g)=\prod_{j=0}^{p-1}f(g\tl t^j)$. Therefore $(\phi_p.f)(s)=\rho_s^p=\chi(s)$ for every $s\in G^{C_p}$. If $g=g_i\tl t^j$ for some $i,j$, then a calculation $(\phi_p.f)(g)=f(g_i)=\chi(g_i)=\chi(g_i\tl t^j)=\chi(g)$, which uses the fact that $\chi$ is a fixed point under the action by $C_p$, completes the proof.

(iv) follows immediately from $H^2_{\c}(\tl)=B^2_{\c}/B^2_c(\tl)$ and parts (i)-(iii).\qed
\begin{Cor}\label{basicisic} There is a $C_p$ isomorphism $H^2_c(\tl)\simeq Z^2_N(\tl)/{\ker\Phi}.$\end{Cor}
\pf Combining Lemmas \ref{hf=adm} and \ref{structureofcoboundaries} with the natural epimorphism $Z^2_N(\tl)\to Z^2_N(\tl)/{\ker\Phi}$ proves the Corollary.\qed
 
We proceed to the main result of the section.
\begin {Prop}\label{hopf2cohomology} Suppose $G$ is a finite abelian $p$-group. If $p$ is odd, or $p=2$ and either $C_2$-action is trivial, or $G$ is an elementary $2$-group, there exists a $C_p$-isomorphism
\begin{equation}\label{structureofcohomology}H^2_c(\tl)\simeq H^2(C_p,\G,\bullet)\times H^2_N(G,\bu{\k}).\end{equation}
\end{Prop}
\pf (1) First we take up the odd case. By the preceeding Corollary we need to decompose $Z^2_N(\tl)/{\ker\Phi}$. We note that for any $p$ there is a group splitting $Z^2(G,\bu{\k})=B^2(G,\bu{\k})\times H^2(G,\bu{\k})$ due to the fact that the group of $1$-cocycles ${\bu{\k}}^{G}$ is injective, and hence so is $B^2(G,\bu{\k})$. We aim at finding a $C_p$-invariant complement to $B^2(G,\bu{\k})$. To this end we recall a well-known isomorphism $\un{a}:H^2(G,\bu{\k})\tilde{\to}\text{Alt}(G)$, see e.g. \cite[\S 2.3]{Y}. There $\text{Alt}\,(G)$ is the group of all bimultiplicative alternating functions 
\begin{equation*}\beta:G\times G\to\bu{\k},\,\beta(ab,c)=\beta(a,c)\beta(b,c),\text{and}\,\beta(a,a)=1\,\text{for all}\,a\in G.\end{equation*}
For future applications we outline the construction of $\un{a}$. Namely, $\un{a}$ is the antisymmetrization mapping sending $z\in Z^2(G,\bu{\k})$ to $\un{a}(z)$ defined by $\un{a}(z)(a,b)=z(a,b)z^{-1}(b,a)$. One can check that $\un{a}$ is bimultiplicative (cf. \cite[(10)]{Y}) and it is immediate that $\un{a}$ is $C_p$-linear. Another verification gives $\text{im}\,\un{a}=\text{Alt}(G)$ and, moreover, $\ker\,\un{a}=B^2(G,\bu{\k})$, see \cite[Thm.2.2]{Y}.  Thus we obtain a $C_p-$isomorphism $H^2(G,\bu{\k})\simeq\text{Alt}(G)$.

For every $\beta=\un{a}(z)$ a simple calculation gives $\un{a}(\beta)=\beta^2$. Since elements of $\text{Alt}(G)$ are bimultiplicative mappings, they have orders dividing the exponent of $G$. Thus $\un{a}(\beta)\ne 1$ for all $\beta\in\text{Alt}(G)$. It follows $B^2(G,\bu{\k})\cap\text{Alt}(G)=\{1\}$. We arrive at a splitting of abelian groups
\begin{equation*}Z^2(G,\bu{\k})=B^2(G,\bu{\k})\times \text{Alt}(G)\end{equation*}
But now both subgroups $B^2(G,\bu{\k})$ and $\text{Alt}(G)$ are $C_p$-invariant hence there holds $Z^2_N(G,\bu{\k})=B^2_N(G,\bu{\k})\times \text{Alt}_N(G)$ which, in view of $\text{Alt}(G)=H^2(G,\bu{\k})$, is the same as
\begin{equation}\label{Z^2adm} Z^2_N(\tl)=B^2_N(\tl)\times H^2_N(G,\bu{\k}).\end{equation} 
Now part (iii) of Lemma \ref{structureofcoboundaries} completes the proof of part (1).

(2) Here we prove the second claim of the Proposition. We decompose $G$ into a product of cyclic groups $\gen{x_i},1\le i\le m$. For every $\alpha\in\text{Alt}(G)$ we define $s_{\alpha}\in Z^2(G,\k^{\bullet})$ via
$$s_{\alpha}(x_i,x_j)=\begin{cases}\alpha(x_i,x_j),&\text{if}\;i\le j\\\phantom{\alpha}1,&\text{else}.\end{cases}$$
The set $S=\{s_{\alpha}|\alpha\in\text{Alt}(G)\}$ is a subgroup. One can see easily that $\un{a}(s_{\alpha})=\alpha$, hence $S$ is isomorphic to $\text{Alt}(G)$ under $\un{a}$. Let us write $G_{(p)}$ for the set of elements of order $p$. We observe that $Z^2_N(\t)=Z^2(\t)_{(2)}$, hence $S_{(2)}\subset Z^2_N(\t)$. For every $z\in Z^2_N(\t),\un{a}(z)\in\text{Alt}(G)_{(2)}$, and therefore $\un{a}(z)=\un{a}(s)$ for some $s\in S_{(2)}$. We have $zs^{-1}\in B^2(G,\k^{\bullet})$, and, as $zs^{-1}$ has order $2$, $zs^{-1}\in B^2_N(\t)$. Thus $Z^2_N(\t)=B^2_N(\t)\times S_{(2)}$ which implies \eqref{structureofcohomology} as $S_{(2)}\simeq\text{Alt}(G)_{(2)}\simeq H^2_N(\t)$.

(3) We prove the last claim of the Proposition. Below $G$ is an elementary $2$-group, and action of $C_2$ is nontrivial. First we establish an intermediate result, namely
\begin{Lem}\label{nonsplitcase} If action $\tl$ is nontrivial, then $Z^2_N(\tl)$ is a nonsplit extension of $\text{Alt}_N(G)$ by $B^2_N(\tl)$.\end{Lem}
\pf This will be carried out in steps.

(i) We aim at finding a basis for $\text{Alt}_N(G)$. We begin by noting that as $\text{Alt}(G)$ has exponent $2$, $\text{Alt}_N(G)$ is the set of all fixed points in $\text{Alt}(G)$. Put $R=\mathbb{Z}_2C_2$. One can see easily that $R$-module $G$ decomposes as
\begin{equation}\label{decomposition}G=R_1\times\cdots\times R_m\times G_0\end{equation}
where  $R_i\simeq R$ as a right $C_2$-module, and $G_0=G^{C_2}$. Denote by $t$ the generator of $C_2$. For each $i$ let $\{x_{2i-1},x_{2i}\}$ be a basis of $R_i$ such that $x_{2i-1}\tl t=x_{2i}$. We also fix a basis $\{x_{2m+1},\ldots,x_n\}$ of $G_0$.

We associate to every subset $\{i,j\}$ the bilinear form $\alpha_{ij}$ by setting
\begin{equation*}\alpha_{ij}(x_i,x_j)=\alpha_{ij}(x_j,x_i)=-1,\,\text{and}\,\alpha_{ij}(x_k,x_l)=1\,\text{for any}\,\{k,l\}\ne\{i,j\}.\end{equation*}
The set $\{\alpha_{ij}\}$ forms a basis of $\text{Alt}(G)$. One can check easily that $t$ acts on basic elements as follows
\begin {equation}\label{actioninAlt} t.\alpha_{ij}=\alpha_{kl}\;\text {if and only if}\,\{x_i,x_j\}\tl t:=\{x_i\tl t,x_j\tl t\}=\{x_k,x_l\}.\end{equation}
Recall the element $\phi_2=1+t\in\mathbb{Z}C_2$. We define forms $\beta_{ij}$ via
\begin{equation}\label{basisAlt}\beta_{ij}=\phi_2.\alpha_{ij}\;\text{if}\;t.\alpha_{ij}\ne\alpha_{ij},\;\text{and}\;\beta_{ij}=\alpha_{ij},\;\text{otherwise}.\end{equation} 
The label $ij$ on $\beta_{ij}$ is not unique as $\beta_{ij}=\beta_{kl}$ whenever $\{x_i,x_j\}\tl t=\{x_k,x_l\}$. Of the two sets $\{i,j\}$ and $\{k,l\}$ labeling $\beta_{ij}$ we agree to use the one with the smallest element, and call such minimal. We claim:
\begin{equation}\label{basisALTN} \text{The elements}\; \{\beta_{ij}\}\; \text{form a basis of}\; \text{Alt}_N(G).\end{equation}
\pf  Suppose $\beta\in\text{Alt}_N(G)$. Say $\beta=\prod\alpha_{ij}^{e_{ij}},\;e_{ij}=0,1$. From $t.\beta=\prod(t.\alpha_{ij})^{e_{ij}}=\beta$ we see that if $\alpha_{ij}$ occurs in $\beta$, i.e. $e_{ij}=1$, then so does $t.\alpha_{ij}$, hence $\beta$ is a product of $\beta_{ij}$.\qed

(ii) We want to show $\text{Alt}_N(G)$ is an epimorphic image of $Z^2_N(\tl)$. 

The restriction $\un{a}^*$ of $\un{a}$ to $Z^2_N(\tl)$ induces a $C_2$-homomorphism $Z^2_N(\tl)\overset{\un{a}^*}\rightarrow\text{Alt}_N(G)$. We have $\ker\,\un{a}^*=B^2(G,\k^{\bullet})\cap Z^2_N(\tl)=B^2_N(\tl)$. First we show $\phi_2.\text{Alt}(G)\subset \text{im}\,\un{a}^*$. For, if $\beta=\phi_2.\alpha$, pick an $s\in Z^2(G,\k^{\bullet})$ with $\un{a}(s)=\alpha$. Then $(t-1).s\in Z^2_N(\tl)$, and $\un{a}((t-1).s)=(t-1).\un{a}(s)=(t-1).\alpha=\phi_2.\alpha=\beta$, as $\alpha^2=1$, which gives the inclusion.

By step (i) and definition \eqref{basisAlt} it remains to show that all fixed points $\alpha_{ij}$ lie in $\text{im}\,\un{a}^*$. By formula \eqref{actioninAlt} $\alpha_{ij}$ is a fixed point if and only if either (a) $\{i,j\}\subset\{2m+1,\ldots,n\}$ or (b) $\{i,j\}=\{2k-1,2k\}$ for some $k,\,1\le k\le m$. Below we find it convenient to write $s_{i,j}$ for $s_{\alpha_{ij}}$.

Consider case (a). We claim $s_{i,j}$ is a fixed point. For, $t.s_{i,j}$ is bimultiplicative, hence is determined by its values at $(x_k,x_l)$. It is immediate that $t.s_{i,j}(x_k,x_l)=s_{i,j}(x_k,x_l)$ for all $(x_k,x_l)$, whence the assertion. Since $s_{i,j}^2=1$ for all $i,j$, $\phi_2.s_{i,j}=1$, hence $s_{i,j}\in Z^2_N(\tl)$. As $\un{a}(s_{i,j})=\alpha_{ij}$, this case is done.

We take up (b). Say $z=s_{2i-1,2i}$ for some $i, 1\le i\le m$. An easy verification gives $\phi_2.z=\alpha_{2i-12i}\ne 1$. Thus $z\notin Z^2_N(\tl)$. To prove (ii) we need to find a coboundary $\delta g_i$ such that $z\delta g_i\in Z^2_N(\tl)$. Since $\un{a}(\alpha_{2i-12i})=1,\alpha_{2i-12i}=\delta f_i$ for some $f_i:G\to\bu{\k}$. Put $G_i$ for the subgroup of $G$ generated by all $x_j,j\ne 2i-1,2i$. We assert that one choice is the function $f_i$ defined by
\begin{equation}\label{specialcocycle} f_i(x_{2i-1}^{j_1}x_{2i}^{j_2}x')=(-1)^{j_1+j_2+j_1j_2}\,\text{for all}\;x'\in G_i\end{equation}
For, on the one hand it is immediate that for any $x',x''\in G_i$ \begin{equation*}\alpha_{2i-12i}(x_{2i-1}^{j_1}x_{2i}^{j_2}x',x_{2i-1}^{k_1}x_{2i}^{k_2}x'')=(-1)^{j_1k_2+j_2k_1}\end{equation*}
On the other hand the definitions of $f_i$ and differential $\delta$ give
\begin{eqnarray*}\lefteqn{\delta  f_i(x_{2i-1}^{j_1}x_{2i}^{j_2}x',x_{2i-1}^{k_1}x_{2i}^{k_2}x'')}\\&&=(-1)^{j_1+j_2+j_1j_2}(-1)^{k_1+k_2+k_1k_2}(-1)^{j_1+k_1+j_2+k_2+(j_1+k_1)(j_2+k_2)}\\&&=(-1)^{j_1k_2+j_2k_1}\nonumber\end{eqnarray*}
Define the function $g_i:G\to\bu{\k}$ by $g_i(x_{2i-1}^{j_1}x_{2i}^{j_2}x')=\iota^{j_1+j_2+j_1j_2}$ where $\iota^2=-1$. One can check easily the equalities $f_i^2=1$ and $t.g_i=g_i,\,g_i^2=f_i$. Hence we have $f_i(\phi_2.g_i)=f_ig_i^2=f_i^2=1$, and then a calculation 
\begin{equation*} \phi_2.(z\delta g_i)=(\phi_2.z)(\phi_2.\delta g_i)=\delta f_i\cdot\delta(\phi_2.g_i)=\delta(f_i(\phi_2.g_i))=1\end{equation*}
completes the proof of (ii).

(iii) Suppose $Z^2_N(\tl)=B^2_N(\tl)\times C$ where $C$ is a $C_2$-invariant subgroup. Then $C$ is mapped isomorphically on $\text{Alt}_N(G)$ under $\un{a}$ and so there is a unique $z\in C$ such that $\un{a}(z)=\alpha_{12}$. Since $\un{a}(s_{1,2})=\alpha_{12},\,z=s_{1,2}\delta g$ for some $g:G\to\bu{\k}$. Further, as $\alpha_{12}$ is a fixed point $\un{a}(t.z)=\alpha_{12}$ as well, hence $t.z=z$. In addition, since $\text{Alt}(G)$ is an elementary $2$-group, 
$1=z^2=(s_{1,2}\delta g)^2=(\delta g)^2=\delta(g^2)$. It follows that $g^2$ is a character of $G$. Moreover, $t.z=z$ is equivalent to $t.s_{1,2}(t.\delta g)=s_{1,2}\delta g$ which in turn gives $s_{1,2}(t.s_{1,2})(t.\delta g)=\delta g$. As $\phi_2.s_{1,2}=\alpha_{12}=\delta f_1$ we have $\delta f_1(t.\delta g)=\delta g$ which implies $\delta f_1=\delta g(t.\delta g)$ on the account of $(\delta g)^2=\delta(g^2)=1$ as $g^2$ is a character. Equivalently we have the equality 
\begin{equation}\label{splitting} f_1=g\cdot(t.g)\cdot\chi\;\text{for some}\;\chi\in\G.\end{equation}
Noting that $f_1$ is defined up to a character of $G$ we can assume that $f_1(x_1)=1=f_1(x_2)\;\text{and}\;f_1(x_1x_2)=-1$. For, $f_1$ is defined as any function satisfying $\delta f_1=\alpha_{12}$. As $\delta(f_1\chi)=\delta f_1$ for any $\chi\in\G$, $f_1$ can be modified by any $\chi$. By \eqref{specialcocycle} $f_1(x_j)=-1=f_1(x_1x_2),j=1,2$ so we can take $\chi$ such that $\chi(x_1)=\chi(x_2)=-1$. The equality \eqref{splitting} implies that for some $\chi\in\G$ there holds  
\begin{align}1&=f_1(x_j)=g(x_1)g(x_2)\chi(x_j),\,j=1,2\tag{*},\, \text{and}\\ 
-1&=f_1(x_1x_2)=g(x_1x_2)^2\chi(x_1x_2)\tag{**}\end{align} 
as $t$ swaps $x_1$ and $x_2$. Since $g^2$ is a character, $g^2(a)=\pm 1$ for every $a\in G$. It follows that $g(x_1)=\iota^m$ and $g(x_2)=\iota^k$ for some $0\le m,k\le 3$. Then equation (*) gives $1=\iota^{m+k}\chi(x_j)$. This equality shows that $\chi(x_1)=\chi_(x_2)$ and $m+k$ is even, because $\chi(a)=\pm 1$ for all $a$. Now (**), and the fact that $g^2$ is a character, gives $-1=g^2(x_1)g^2(x_2)\chi(x_1)\chi(x_2)=\iota^{2(m+k)}\iota^{-2(m+k)}=1$, a contradiction. This completes the proof of the Lemma.\qed

Finally we prove (3). Let $G$ be a group with a decomposition \eqref{decomposition}. Set $C$ to be the subgroup of $Z^2_N(\tl)$ generated by the set $B=B'\cup B''\cup B'''$ where
\begin{align*}B'&=\{\phi_2.s_{i,j}|\alpha_{ij}\,\text{is not a fixed point, and}\,\{ij\}\,\text{is minimal}\} \\B''&=\{s_{i,j}|i<j\,\text{and}\,\{i,j\}\subset\{2m+1,\ldots,n\}\}\\B'''&=\{s_{2i-1,2i}\delta g_i|i=1,\ldots,m\}.\end{align*} 
There $\delta g_i$ is chosen as in the proof of the case (ii) of Lemma \ref{nonsplitcase}. Passing on to $Z^2_N(\tl)/\ker\Phi$ we denote by $\ov{B^2_N(\tl)}$ and $\ov{C}$ images of these subgroups in $Z^2_N(\tl)/\ker\Phi$. Pick a $v\in B$. If $v\in B'\cup B''$ then $v^2=1$ because the corresponding $s_{i,j}$ has order $2$. For $v=s_{2i-1,2i}\delta g_i$, $v^2=\delta g_i^2=\delta f_i$. We know $t.f_i=f_i$ and $f_i^2=1$ and therefore $\phi_2.f_i=1$, whence $\delta f_i\in\ker\Phi$ by definition \eqref{coboundarymap}. It follows that ${\ov{v}}^2=1$ for all $\ov{v}\in\ov{B}$. Furthermore, by Lemma \ref{nonsplitcase} the mapping $\un{a}$ sends $\ov{B}$ to the basis \eqref{basisALTN} of $\text{Alt}_N(G)$. Therefore $\ov{C}$ is isomorphic to $\text{Alt}_N(G)$ at least as an abelian group and forms a complement to $\ov{B^2_N(\tl)}$ in $Z^2_N(\tl)/\ker\Phi$. Since $\text{Alt}_N(G)$ consists of fixed points the proof will be completed if we show the same for $\ov{C}$. The fact that $B'\cup B''$ consists of fixed points follows from $t\phi_2=\phi_2$ and the case (i) of Lemma \ref{nonsplitcase}. For an $s_{2i-1,2i}\delta g_i$, the equality $\phi_2.s_{2i-1,2i}=s_{2i-1,2i}(t.s_{2i-1,2i})=\alpha_{2i-12i}=\delta f_i$ gives $t.s_{2i-1,2i}=s_{2i-1,2i}\delta f_i$. Since $\delta f_i\in\ker\Phi$ and $t.\delta g_i=\delta g_i$ we see that $\overline{s_{2i-1,2i}\delta g_i}$ is a fixed point in $Z^2_N(\tl)/\ker\Phi$ which completes the proof.\qed 
 
\section{The Isomorphism Theorems}\label{Main}
We begin with a general observation. Let $H$ be be an extension of type (A). The mapping $\pi$ induces a $\k F$-comodule structure $\rho_{\pi}$ on $H$ via
\begin{equation}\label{comodulealg}\rho_{\pi}:H\to H\o\k F,\,\rho_{\pi}(h)=h_1\o\pi(h_2).\end{equation}
$H$ becomes an $F$-graded algebra with the graded components $H_f=\{h\in H|\rho_{\pi}(h)=h\o f\}$. Let $\chi:\k F\to H$ be a section of $\k F$ in $H$. By definition $\chi$ is a convolution invertible $\k F$-comodule mapping, that is
\begin{equation}\label{Section}\rho_{\pi}(\chi(f))=\chi(f)\o f,\,\text{for every}\,f\in F\end{equation}
Set $\ov{f}=\chi(f)$. The next lemma is similar to \cite[3.4]{M} or \cite[7.3.4]{Mo}.
\begin{Lem}\label{gradedcomponents} For every $f\in F$ there holds $H_f=\D{G}\ov{f}$\end{Lem}
\pf By definition of components $H_1=H^{\text{co}\pi}$ which equals to $\D{G}$ by the definition of extension. By the equation \eqref{Section} $\rho_{\pi}(\ov{f})=\ov{f}\o f$, hence $\D{G}\ov{f}\subset H_f$. Since the containment holds for all $f$, the equalities 
\begin{equation*}H=\oplus_{f\in F}H_f=\oplus_{f\in F}\D{G}\ov{f}\end{equation*}
force the equalities $H_f=\D{G}\ov{f}$ for all $f\in F$.\qed
\begin{Def}Given two $F$-graded algebras $H=\oplus H_f$ and $H'=\oplus H'_f$ and an automorphism $\alpha:F\to F$ we say that a linear mapping $\psi:H\to H'$ is an $\alpha$-graded morphism if $\psi(H_f)=H'_{\alpha(f)}$ for all $f\in F$.\end{Def}
\begin{Lem}\label{gradedpsi} Suppose $H$ and $H'$ are two extensions of $\k F$ by $\D{G}$ and $\psi:H\to H'$ a Hopf isomorphism sending $\D{G}$ to $\D{G}$. Then $\psi$ is an $\alpha$- graded mapping for some $\alpha$.\end{Lem}
\pf Suppose $H$ and $H'$ are given by sequences
\begin{equation*}\k^G\overset\iota\rightarrowtail H\overset\pi\twoheadrightarrow \k F,\,\text{and}\quad\k^G\overset{\iota'}\rightarrowtail H'\overset{\pi'}\twoheadrightarrow \k F\end{equation*}
By definition of extension $\text{Ker}\,{\pi}=H(\D{G})^+$ and likewise $\text{Ker}\,{\pi'}=H'(\D{G})^+$. By assumption $\psi(\D{G})=\D{G}$, hence $\psi$ induces a Hopf isomorphism $\alpha:H/H(\D{G})^+\to H'/H'(\D{G})^+$. Replacing $H/H(\D{G})^+$ and $ H'/H'(\D{G})^+$ by $\k F$ we can treat $\alpha$ as a Hopf isomorphism $\alpha:\k F\to\k F$. $\alpha$ is in fact an automorphism of $F$. We arrive at a commutative diagram
\begin{equation*}\begin{CD}
{\k}^G @>\iota>> H @>\pi>>\k F\\
@V{\psi}VV @V{\psi}VV  @V{\alpha}VV \\      
{\k}^G @>\iota'>> H' @>\pi'>>\k F
\end{CD}\end{equation*}
Since $\psi$ is a coalgebra mapping for every $f\in F$ we have 
\begin{align*}\Delta_{H'}(\psi(\ov{f}))&=(\psi\o\psi)\Delta_H(\ov{f})=\psi((\ov{f})_1)\o\psi((\ov{f})_2),\,\text{hence}\\
\rho_{\pi'}(\psi(\ov{f}))&=\psi((\ov{f})_1)\o \pi'\psi((\ov{f})_2)=\psi((\ov{f})_1)\o\alpha\pi((\ov{f})_2)\end{align*}
On the other hand, applying $\psi\o\alpha$ to the equality 

\noindent$\rho_{\pi}(\ov{f})=(\ov{f})_1\o\pi((\ov{f})_2)=\ov{f}\o f$ gives
\begin{equation*}\psi((\ov{f})_1)\o\alpha\pi((\ov{f})_2)=\psi(\ov{f})\o\alpha(f)\end{equation*}
whence we deduce $\rho_{\pi'}(\psi(\ov{f}))=\psi(\ov{f})\o\alpha(f)$. Thus $\psi(\ov{f})\in H'_{\alpha(f)}$ which shows the inclusion 
\begin{equation*}\psi(H_f)=\psi(\D{G}\ov{f})=\D{G}\psi(\ov{f})\subseteq H'_{\alpha(f)}=\D{G}\ov{\alpha(f)}\end{equation*}
Since both sides of the above inclusion have equal dimensions, the proof is complete.\qed

From this point on $F=C_p$. Let $\tl$ and $\tl'$ be two actions of $C_p$ on $G$. We denote $(G,\tl)$ and $(G,\tl')$ the corresponding $C_p$-modules and we use the notation `$\bullet$' and `$\circ$' for the actions of $C_p$ on $\D{G}$ corresponding by \eqref{dualFaction} to $\tl$ and $\tl'$, respectively. We let $I(\tl,\tl')$ denote the set of all automorphisms of $G$ intertwining actions $\tl$ and $\tl'$, that is automorphisms $\lambda:G\to G$ satisfying
\begin{equation}\label{modulemap}(a\tl x)\lambda=a\lambda\tl' x,\,a\in G,x\in C_p\end{equation}
We make every $\lambda$ act on  functions $\tau:C_p\times G^2\to\D{C_p}$ by 
\begin{equation*}(\tau.\lambda)(x,a,b)=\tau(x,a\lambda^{-1},b\lambda^{-1}).\end{equation*} 

\begin{Lem}\label{actiononcocycles} {\rm(i)} The group $Z^2(G,\bu{(\D{C_p})})$ is invariant under the action induced by any automorphism of $G$,

{\rm(ii)} A $C_p$-isomorphism $\lambda:(G,\tl)\to (G,\tl')$ induces $C_p$-isomorphisms between the groups $Z^2_c(\tl),B^2_c(\tl),H^2_c(\tl)$ and $Z^2_c(\tl'),B^2_c(\tl'),H^2_c(\tl')$, respectively.
\end{Lem}
\pf (i) is immediate.

(ii) We must check condition \eqref{hopfcocycle} for $\tau.\lambda$ and $\mathbb{Z}C_p$-linearity of the induced map. First we note $\lambda^{-1}$ is a $C_p$- isomorphism between $(G,\tl')$ and $(G,\tl)$, as one can check readily. For, set $b=a\lambda$ in \eqref{modulemap}. Then we have $(b\lambda^{-1}\tl x)\lambda=(b\tl'x)$ hence $b\lambda^{-1}\tl x=(b\tl' x)\lambda^{-1}$

Next we verify \eqref{hopfcocycle} and $C_p$-linearity in a single calculation
\begin{align*} &(\tau.\lambda)(xy)(a,b)=\tau(xy,a\lambda^{-1},b\lambda^{-1})\\
&=\tau(x,a\lambda^{-1},b\lambda^{-1})(x\bullet\tau(y,a\lambda^{-1},b\lambda^{-1}))\\
&=\tau(x,a\lambda^{-1},b\lambda^{-1})\tau(y,a\lambda^{-1}\tl x,b\lambda^{-1}\tl x)\\&=\tau(x,a\lambda^{-1},b\lambda^{-1})\tau(y,(a\tl' x)\lambda^{-1},(b\tl' x)\lambda^{-1})\\
&=(\tau.\lambda)(x)(x\circ(\tau.\lambda)(y))(a,b).
\end{align*}                                                                                                                                                 In the case of $B^2_c(\tl)$, first one checks the equality
\begin{equation*}(\delta_{G}\eta).\lambda=\delta_{G}(\eta.\lambda)\,\text{for any}\,\eta:C_p\times G\to\bu{\k}.\end{equation*}                                                        
It remains to verify the condition \eqref{hopf1cocycle} for $\eta.\lambda$. That is done similarly to the calculation in (ii).\qed

Let $(G,\tl)$ be a $C_p$-module. We denote by $\mathbb{A}(\tl)$ the group of $C_p$-automorphisms of $(G,\tl)$. By the above Lemma $Z^2_c(\tl)$ is an $\mathbb{A}(\tl)$-module. Symmetrically, we introduce the group $A_p=\mathrm{Aut}(C_p)$ of automorphisms of $C_p$. We define an action of $A_p$ on $\mathrm{Map}(C_p\times G^2,\bu{\k})$ via
\begin{equation*}\tau.\alpha(x,a,b)=\tau(\alpha(x),a,b)\end{equation*}
We want to know the effect of this action on $Z^2_c(\tl)$. Let $(G,\tl)$ be a $C_p$-module. For $\alpha\in A_p$ we define a $C_p$-module $(G,\tl^{\alpha})$ via
\begin{equation*}a\tl^{\alpha}x=a\tl\alpha(x),\;a\in G,\,x\in C_p\end{equation*}
Similarly, an action `$\bullet$' of $C_p$ on $\k^G$ can be twisted by $\alpha$ into `$\bullet^{\alpha}$' by 
\begin{equation*}x\bullet^{\alpha}r=\alpha(x)\bullet r,\;r\in\k^G\end{equation*}  
One can see easily that if $\bullet$ and $\tl$ correspond to each other by \eqref{dualFaction}, then so do $\bullet^{\alpha}$ and $\tl^{\alpha}$. 
\begin{Lem}\label{secondaction} {\rm(i)} If $\lambda\in I(\tl,\tl')$, then $\lambda\in I(\tl^{\alpha},\tl'^{\alpha})$ for every $\alpha\in A_p$,
 
\rm{(ii)} The mapping $\tau\mapsto\tau.\alpha$ induces an $\mathbb{A}(\tl)$-isomorphism between $Z^2_c(\tl),B^2_c(\tl),H^2_c(\tl)$ and $Z^2_c(\tl^{\alpha}),B^2_c(\tl^{\alpha}),H^2_c(\tl^{\alpha})$, respectively for every $\alpha\in A_p$ .\end{Lem}
\pf {\rm(i)} For every $a\in G,x\in C_p$ we have
\begin{equation*}(a\tl^{\alpha})\lambda=(a\tl\alpha(x))\lambda=a\lambda\tl'\alpha(x)=a\lambda\tl'^{\alpha}x\end{equation*}

{\rm(ii)} First we note that $\mathbb{A}(\tl)$ can be identified with $\mathbb{A}(\tl^{\alpha})$ for any $\alpha$ by the folllowing calculation
\begin{equation*}(g\tl^{\alpha}x)\phi=(g\tl\alpha(x))\phi=(g\phi)\tl\alpha(x)=g\phi\tl^{\alpha}x\,\text{for every}\,\phi\in\mathbb{A}(\tl).\end{equation*}
Thus we will treat every $Z^2_c(\tl^{\alpha})$ as an $\mathbb{A}(\tl)$- module. Our next step is to show that for every $\tau\in Z^2_c(\tl)$, $\tau.\alpha$ lies in $Z^2_c(\tl^{\alpha})$. This boils down to checking \eqref{hopfcocycle} for $\tau.\alpha$ with the $\tl^{\alpha}$-action: 
\begin{align*}(\tau.\alpha)(xy)&=\tau(\alpha(x)\alpha(y))=\tau(\alpha(x))(\alpha(x)\bullet\tau(\alpha(y))\\&=\tau(\alpha(x))(x\bullet^{\alpha}\tau(\alpha(y))=(\tau.\alpha)(x)(x\bullet^{\alpha}(\tau.\alpha)(y)).\end{align*}
As for $\mathbb{A}(\tl)$-linearity, for every $\phi\in\mathbb{A}(\tl)$, we have
\begin{align*}((\tau.\alpha).\phi)(x,a,b)&=(\tau.\alpha)(x,a\phi^{-1},b\phi^{-1})=\tau(\alpha(x),a\phi^{-1},b\phi^{-1})\\&=(\tau.\phi)(\alpha(x),a,b)=((\tau.\phi).\alpha)(x,a,b).\end{align*}\qed

We need several short remarks.
\begin{Lem}\label{correction} Suppose $\tau$ is a $2$-cocycle. Assume $r\in\bu{(\D{G})}$ is such that $\phi_p.r=\epsilon$. Set $r_i=\phi_i.r,\,1\le i\le p$. Define a $1$-cocycle $\zeta:G\to\bu{(\D{C_p})}$ by $\zeta(t^i)=r_i$ and a $2$-cocycle $\tau'=\tau(\delta_{G}\zeta)$. Then the mapping
\begin{equation*}\iota:H(\tau,\tl)\to H(\tau',\tl),\,\iota(p_at^i)=p_ar_it^i,\,a\in G,\,1\le i\le p\end{equation*}
is an equivalence of extensions.\end{Lem}
\pf We need to show $\delta_{G}\zeta\in B^2_c$ which means that $\zeta$ satisfies \eqref{hopf1cocycle}. The argument of Lemma \ref{hf=adm} used to derive \eqref{hopfcocycle} from the condition \eqref{Adm} works verbatim for $\zeta$.\qed

\begin{Lem}\label{cocommutative} $H(\tau,\tl)$ is cocommutative  iff $\tau$ lies in $H^2_{\c}(\tl)$.\end{Lem}
\pf $H^*(\tau,\tl)$ is commutative iff $\ov{a}\ov{b}=\ov{b}\ov{a}$ which is equivalent to $\tau(a,b)=\tau(b,a)$. The latter implies that $\tau(t):G\times G\to\bu{\k}$ is a symmetric $2$-cocycle, hence a coboundary, that is an element of $B^2_N$. A reference to Lemma \ref{structureofcoboundaries}(i) completes the proof.\qed

Unless stated otherwise, $H(\tau,\tl)$ is a noncocommutative Hopf algebra. We pick another algebra $H(\tau',\tl')$ isomorphic to $H(\tau,\tl)$ via $\psi:H(\tau,\tl)\to H(\tau',\tl')$. The next observation is noted in \cite[p. 802]{Mas1}. 
\begin{Lem}\label{Gstability} Mapping $\psi$ induces an Hopf automorphism of $\D{G}$.\end{Lem}\qed

Let $G$ be a finite group and $\mathrm{Aut_{Hf}}(\D{G})$ be the group of Hopf automorphisms of $\D{G}$. For $\phi\in\mathrm{Aut_{Hf}}(\D{G})$ we denote by $\phi^*$ the mapping of $G$ induced by $\phi$ via
\begin{equation}\label{dualphi}(a\phi^*)(f):=f(a\phi^*)=\phi(f)(a),\,f\in\D{G}.\end{equation} 
\begin{Lem}\label{dualofphi} Let $G$ be a finite abelian group. The mapping $\phi\mapsto\phi^*$ is an isomorphism between $\mathrm{Aut_{Hf}}(\D{G})$ and $\mathrm{Aut}(G)$. $\phi$ is a $C_p$-isomorphism $(\D{G},\bullet)\to (\D{G},\circ)$ if and only if $\phi^*$ is a $C_p$-isomorphism $(G,\tl')\to (G,\tl)$.\end{Lem}
\pf In general $\phi^*$ is a permutation of the set $G$. When $G$ is abelian and $\k$ contains a $|G|$th root of $1$, we have $\D{G}=\k\G$. Then $\phi(\G)=\G$, as $\phi$ preserves grouplikes. It follows from a straightforward calculation that $\phi^*$ is a group automorphism and $\phi\mapsto\phi^*$ is an isomorphism.

We proceed to formulation of isomorphism theorems. We need several preliminary remarks. First off, let $\tl$ be a $C_p$-action on $G$. We denote by $[\tl]$ the class of $C_p$-actions $\tl'$ isomorphic to $\tl^{\alpha}$ for some $\alpha\in A_p$, that is such that $I(\tl',\tl^{\alpha})$ is nonempty. We let $\mathrm{Ext}_{[\tl]}(\k C_p,\k^G)$ stand for all equivalence classes of extensions whose $C_p$-action on $G$ lies in $[\tl]$.

In the second place we construct groups $\g(\tl)$ that control isomorphism types of extensions. For the trivial action we set $\g(\t)=\mathrm{Aut}(G)\times A_p$. Else, we observe that by Lemma \ref{secondaction}(i) if $\lambda\in I(\tl,\tl^{\alpha}),\mu\in I(\tl,\tl^{\beta})$, then $\lambda\mu\in I(\tl,\tl^{\alpha\beta})$. Therefore the set of all $\alpha\in A_p$ such that $I(\tl,\tl^{\alpha})\ne\emptyset$ is a subgroup of $A_p$ denoted by $C(\tl)$. Let us select an element $\lambda_{\alpha}\in I(\tl,\tl^{\alpha})$ for every $\alpha\in C(\tl)$. We define $\g(\tl)$ as the subgroup of $\mathrm{Aut}(G)$ generated by $\A$ and the elements $\lambda_{\alpha},\alpha\in C(\tl)$.
\begin{Prop}\label{keygroup} If action $\tl$ is nontrivial, then $\g(\tl)$ is a crossed product of $\A$ with $C(\tl)$.\end{Prop}
\pf It is evident that $\lambda\A\lambda^{-1}=\A$ for every $\lambda\in I(\tl,\tl^{\alpha})$. In addition, for every $\lambda,\mu\in I(\tl,\tl^{\alpha}),\,\lambda^{-1}\mu\in\A$. Thus we have $I(\tl,\tl^{\alpha})=\A\lambda_{\alpha}$. It follows that $\lambda_{\alpha}\cdot\lambda_{\beta}=\phi(\alpha,\beta)\lambda_{\alpha\beta}$ for some $\phi(\alpha,\beta)\in\A$. It remains to show that the kernel of $\pi:\g(\tl)\to C(\tl),\pi(\phi\lambda_{\alpha})=\alpha$ equals $\A$. Pick $\alpha:x\to x^k,k\ne 1$. Clearly $\lambda\in I(\tl,\tl^{\alpha})$ iff $t\lambda=\lambda t^k$ where we treat $t\in C_p$ as automorphism of $G$. Since elements of $\A$ commute with $t$, $I(\tl,\tl^{\alpha})\cap\A=\emptyset$.\qed

Our next goal is to define a $\g(\tl)$-module structure on $H^2_c(\tl)$. For every $\lambda\in I(\tl,\tl^{\alpha})$ Lemmas \ref{actiononcocycles}(ii),\ \ref{secondaction}(ii) show that the mapping
\begin{equation}\label{mainautomorphisms}\omega_{\lambda,\alpha}: \tau\mapsto\tau.\lambda\alpha^{-1},\tau\in H^2_c(\tl)\end{equation}
is an automorphism of $H^2_c(\tl)$. For $\lambda=\lambda_{\alpha}$ we write $\omega_{\alpha}=\omega_{\lambda,\alpha}$. $\A$ also acts on $H^2_c(\tl)$, and we denote by $\ov{\phi}$ the automorphism of $H^2_c(\tl)$ induced by $\phi\in\A$.
\begin{Lem}\label{gactiononcohomology} The mapping $\phi\lambda_{\alpha}\mapsto\ov{\phi}\omega_{\alpha},\phi\in\A,\alpha\in C(\tl)$ defines $\g(\tl)$-module structure on $H^2_c(\tl)$.\end{Lem}
\pf $H^2_c(\tl)$ is a subquotient of $Z^2(G,\bu{(\k^{C_p})})$, and actions of $\A$ and $\omega_{\alpha}$ on $H^2_c(\tl)$ are induced from their actions on $Z^2(G,\bu{(\k^{C_p})})$. By Lemma \ref{actiononcocycles}(i) $Z^2(G,\bu{(\k^{C_p})})$ is an $\mathrm{Aut}(G)$-module, hence it is a $\g(\tl)$-module as well. On the other hand, it is elementary to check that every $\lambda\in\mathrm{Aut}(G)$ commutes with every $\beta\in A_p$ as mappings of $Z^2(G,(\k^{C_p})^{\bullet})$. It follows that the equalities $\omega_{\alpha}\omega_{\beta}=\ov{\phi(\alpha,\beta)}\omega_{\alpha\beta}$ and $\omega_{\alpha}\ov{\phi}\omega_{\alpha}^{-1}=\lambda_{\alpha}\ov{\phi}\lambda_{\alpha}^{-1}$ hold in $\mathrm{Aut}(Z^2(G,\bu{(\k^{C_p})})$. This shows that the mapping of the Lemma is a homomorphism, as needed.\qed
\begin{Thm}\label{mainthm} {\rm (I)}. Noncocommutative extensions $H(\tau,\tl)$\\ and $H(\tau',\tl')$ are isomorphic if and only if 
\begin{enumerate}
\item[(i)] There exist $\alpha\in A_p$ and $C_p$-isomorphism $\lambda:(G,\tl)\to(G,\tl'^{\alpha})$ such that 
\item[(ii)] $\tau'=\tau.(\lambda\alpha^{-1})$ in $H^2_c(\tl')$.
\end{enumerate}
{\rm(II)}. There is a bijection between the orbits of $\g(\tl)$ in $H^2_c(\tl)$ not contained in $H^2_{\c}(\tl)$ and the isomorphism classes of noncocommutative extensions in $\mathrm{Ext}_{[\tl]}(\k C_p,\k^G)$.
\end{Thm}
\pf (I). In one direction, suppose $\psi:H(\tau,\tl)\to H(\tau',\tl')$ is an isomorphism. By Lemma \ref{Gstability} $\psi$ induces  an automorphism $\phi:\D{G}\to\D{G}$, and from  Lemma \ref{gradedpsi} we have the equality $\psi(t)=rt^k$ for some $k$ and $r\in\D{G}$. The equality $\psi(t^p)=1$ implies $(rt^k)^p=\phi_p(t^k)\circ r=1$ and, as $\phi_p(t^k)=\phi_p(t)$, we have $\phi_p\circ r=1$. This shows $r\in\bu{(\D{G})}$. Let $\alpha:x\mapsto x^k, x\in C_p$ be this automorphism of $C_p$, and set $\phi=\psi|_{\k^G}$. Then the calculation
\begin{equation*}\phi(t\bullet f)=\psi(tft^{-1})=r\alpha(t)\phi(f)\alpha(t)^{-1}r^{-1}=\alpha(t)\circ\phi(f),\,f\in\D{G}\end{equation*}
shows $\phi:(\D{G},\bullet)\to(\D{G},\circ^{\alpha})$ is a $C_p$-isomorphism. It follows by Lemma \ref{dualofphi} that $(G,{\tl'}^{\alpha})$ is isomorphic to $(G,\tl)$ under $\phi^*$, hence $\lambda=(\phi^*)^{-1}:(G,\tl)\to (G,{\tl'}^{\alpha})$ is a required isomorhism.

It remains to establish the second condition of the theorem. Set $s=\phi^{-1}(r)$ and observe that, as $\phi^{-1}$ is a $C_p$-mapping, $\phi_p\bullet s=1$.
For $\phi_p\circ r=\phi_p\circ^{\alpha}r=1$. Hence $\phi^{-1}(\phi_p\circ^{\alpha}r)=\phi_p\bullet\phi^{-1}(r)=\phi_p\bullet s=1$, and
therefore by Lemma \ref{correction} there is an equivalence $\iota: H(\tau,\tl)\to H(\widetilde{\tau},\tl)$ with $\iota(t)=st$. 

By construction $\iota$ is an algebra map with $\iota(s)=s$ for all $s\in\k^G$. Hence $t=\iota(st)=s\iota^{-1}(t)$ whence $\iota^{-1}(t)=s^{-1}t$. Thus we have $(\psi\iota^{-1})(t)=t^k$ by the choice of $s$. It follows we can assume $\psi(t)=t^k$ hence $\psi(x)=x^k$ for all $x\in C_p$. 

Abbreviating $H(\tau,\tl),H(\tau',\tl')$ to $H,H'$, respectively, we take up the identity. 
\begin{equation*}\Delta_{H'}(\psi(x))=(\psi\o\psi)\Delta_H(x), x\in C_p,\end{equation*}
expressing comultiplicativity of $\psi$ on elements of $C_p$. By \eqref{DeltaH} this translates into
\begin{equation}\label{identity}\sum_{a,b}\tau'(x^k,a,b)p_ax^k\o p_bx^k=\sum_{c,d}\tau(x,c,d)\phi(p_c)x^k\o\phi(p_d)x^k.\end{equation}
Next we connect $\phi(p_b)$ to the action of $\phi^*$. The argument used to prove \eqref{.vstriangleleft} yields
\begin{equation}\label{actionofphi}\phi(p_b)=p_{b(\phi^*)^{-1}}.\end{equation}

For, since $\phi$ is an algebra map, $\phi(p_b)=p_c$ where $c$ is such that $\phi(p_b)(c)=1$. By definition of action $\phi^*$, $\phi(p_b)(c)=(c\phi^*)(p_b)=p_b(c\phi^*)$, hence $c\phi^*=b$, whence $c=b(\phi^*)^{-1}$.

Switching summation symbols $c,d$ to $l=c(\phi^*)^{-1}$ and $m=d(\phi^*)^{-1}$, the right-hand side of \eqref{identity} takes on the form
\begin{equation*} \sum_{l,m}\tau(x,l\phi^*,m\phi^*)p_lx^k\o p_mx^k\end{equation*}
Thus $\psi$ is comultiplicative on $C_p$ iff 
\begin{equation}\label{comultiplicative}\tau'(\alpha(x),a,b)=\tau(x,a\phi^*,b\phi^*)=\tau(\phi^*)^{-1}(x,a,b)=\tau.\lambda(x,a,b).\end{equation}
Applying $\alpha^{-1}$ to the last displayed equation we arrive at
\begin{equation}\label{comultiplicative'}\tau'(x,a,b)=\tau.\lambda\alpha^{-1}(x,a,b).\end{equation}
as needed. 

Conversely, let us assume hypotheses of part (I). Using Lemma \ref{dualofphi} we infer that $\lambda^{-1}$ induces a Hopf $C_p$-isomorphism $\phi:(\D{G},\bullet)\to(\D{G},\circ^{\alpha})$. We define 
\begin{equation*}\psi:H(\tau,\tl)\to H(\tau',\tl')\,\text{via}\,\psi(fx)=\phi(f)\alpha(x),f\in\D{G},x\in C_p.\end{equation*}
A routine verification using $\phi(x\bullet f)=\alpha(x)\circ\phi(f)$ shows $\psi$ is an algebra mapping.
\begin{align*}&\psi((fx)(f'x'))=\psi(f(x\bullet f')xx')=\phi(f)\phi(x\bullet f')\alpha(x)\alpha(x')\\
&=\phi(f)\alpha(x)\circ\phi(f')\alpha(x)\alpha(x')=\phi(f)\alpha(x)\phi(f')\alpha^{-1}(x)\alpha(x)\alpha(x')\\
&=(\phi(f)\alpha(x))(\phi(f')\alpha(x')=\psi(fx)\psi(f'x').\end{align*}
To see comultiplicativity of $\psi$ we need to verify 
\begin{equation}\label{comult}\Delta_{H'}(\psi(fx))=(\psi\otimes\psi)\Delta_H(fx).\end{equation}
By multiplicativity of $\Delta_{H'},\psi,\Delta_H$ it suffices to check \eqref{comult} for $\phi$ and for every $\psi(x)$. Now the first case holds as $\phi$ is a coalgebra mappping, and the second follows from $\tau'=\tau.\lambda\alpha^{-1}$ by calculations \eqref{identity} and \eqref{comultiplicative'}.

(II). Pick an algebra $H(\tau',\tl')$ in $\mathrm{Ext}_{[\tl]}(\k C_p,\k^G)$. Let us define the set $\cC=\cC(\tau',\tl')$ by the formula
\begin{equation*}\cC(\tau',\tl')=\{(\tau'',\tl'')|H(\tau'',\tl'')\simeq H(\tau',\tl')\}.\end{equation*}
Clearly the family of sets $\{\cC\}$ is identical to the set of isoclasses of extensions in $\mathrm{Ext}_{[\tl]}(\k C_p,\k^G)$. We look at the intersection $\cC\cap H^2_c(\tl)$ as $\cC$ runs over $\{\cC\}$. First, we claim that $\cC(\tau',\tl')\cap H^2_c(\tl)\ne\emptyset$ for every $(\tau',\tl')$. To this end we note that as $\tl'\in[\tl]$ there exists $\mu:(G,\tl')\to(G,\tl^{\alpha})$, and then setting $\tau=\tau'.\mu\alpha^{-1}$ we have  $(\tau,\tl)\in\cC(\tau',\tl')$ by part (I). Next we show the equality $\cC(\tau',\tl')\cap H^2_c(\tl)=\tau\g(\tl)$. For, by definition
$(\sigma,\tl)\in\cC(\tau',\tl')$ iff $H(\sigma,\tl)\simeq H(\tau,\tl)$ which by part (I) implies $\sigma=\tau.\omega_{\lambda,\alpha}$ for some $\alpha\in A_p$ and $\lambda:(G,\tl)\to(G,\tl^{\alpha})$. It follows that the mapping 
\begin{equation*}\cC\to\cC\cap H^2_c(\tl),\cC\in\{\cC\}\end{equation*}
is an injection from the set of isoclasses of noncocommutative extensions in $\mathrm{Ext}_{[\tl]}(\k C_p,\k^G)$ to the set of orbits of $\g(\tl)$ in $H^2_c(\tl)$ not contained in $H^2_{\c}(\tl)$. This mapping is also a surjection as for every $\tau\in H^2_c(\tl)$, $\cC(\tau,\tl)\cap H^2_c(\tl)=\tau\g(\tl)$.\qed
\begin{Cor}\label{sizeoforbit} For every $\tau\in H^2_c(\tl)$ the cardinality of the orbit $\tau\g(\tl)$ satisfies
\begin{equation*}|\tau\A|\le|\tau\g(\tl)|\le|\mathbb{C}(\tl)||\tau\A|.\end{equation*}
\end{Cor}
\pf Let $\X/\A$ be the set of $\A$-orbits in $\X$. Since $\A$ is normal in $\g(\tl)$, there is an induced action of $C(\tl)=\g(\tl)/\A$ on $\X/\A$. An orbit of $\g(\tl)$ is union of points of some orbit of $C(\tl)$. The latter has size $\le|C(\tl)|$ whence the Corollary.\qed

The second isomorphism theorem concerns cocommutative extensions in $\mathrm{Ext}_{[\tl]}(\k C_p,\k^G)$ under a stricter condition on $G$, namely we assume $G$ to be an elementary $p$-group.  
\begin{Thm}\label{cocommutativeclasses} Let $G$ be a finite elementary $p$-group. Then there is a bijection between the set of orbits of $\A$ in  $H^2_{\c}(\tl)$ and the set of isoclasses of cocommutative Hopf algebras in $\mathrm{Ext}_{[\tl]}(\k C_p,\D{G})$.\end{Thm}
\pf By Lemma \ref{cocommutative} $\tau=\delta\eta\in H^2_{\c}(\tl)$. A proof of the Theorem comes down to the statement
\begin{equation}\label{isomorphismcocommutative} H(\delta\eta,\tl)\simeq H(\delta\zeta,\tl')\;\text{iff}\;\delta\zeta=(\delta\eta).(\phi^*)^{-1}
\end{equation}
for a $C_p$-isomorphism $\phi:(\G,\bullet)\to (\G,\circ)$. This makes sense as $\phi:(\k^G,\bullet)\to(\k^G,\circ)$ restricts to $\phi:(\G,\bullet)\to (\G,\circ)$ by Lemma \ref{dualofphi} and $\G$ is $C_p$-stable by Lemma \ref{cocentralext}(2).  The proof of \eqref{isomorphismcocommutative} will be based on several intermediate results.

By general principles a cocommutative Hopf algebra $H(\tau,\tl)$ is a Hopf group algebra of some group $L$. Our first step is to identify that group.
\begin{Lem}\label{cocommutativegroup} A cocommutative extension $H(\tau,\tl)$ is isomorphic as a Hopf algebra to a group algebra $\k L$ with $L\in\mathrm{Opext}(C_p,\G,\bullet)$.
\end{Lem}
\pf Let $H\leadsto G(H)$ be the functor of taking the group of grouplikes of $H$. Applying $G(\cdot)$ to an extension $\k^G\rightarrowtail \k L \twoheadrightarrow \k C_p\in\mathrm{Ext}(\k C_p,\k^G,\tl)$ yields an extension $\G\rightarrowtail L\twoheadrightarrow C_p\in \mathrm{Opext}(C_p,\G,\bullet)$.

On the other hand, by Lemma~\ref{structureofcoboundaries} we have a group isomorphism
\begin{equation*}H^2_{\c}(\tl)/B^2_c(\tl)\simeq {\G}^{C_p}/N(\G)\end{equation*}
under the mapping $\delta\eta\mapsto \Phi\Theta(\delta\eta)=\phi_p\bullet\eta(t)N(\G)$. We will write $L(\ov{\chi})$ when the cohomology class of $L$ is $\ov{\chi}:=\chi N(\G),\,\chi\in {\G}^{C_p}$. We want to construct an explicit isomorphism $H(\delta\eta,\tl)\simeq\k L(\Phi\Theta(\delta\eta))$.

We will use the notation $\chi_{\eta}=\phi_p\bullet\eta(t)$. We observe that $\ov{\chi_{\eta}}=\ov{\epsilon}$ means $\delta\eta$ is a Hopf coboundary. In this case 

\noindent$H(\delta\eta,\tl)\simeq H(\epsilon\otimes\epsilon,\tl)$, where $\epsilon\otimes\epsilon$ is the trivial $2$-cocycle. By \eqref{DeltaH} we have in $H(\epsilon\otimes\epsilon,\tl)$
\begin{equation*}\Delta_H(t)=\sum_{a,b\in G}p_at\otimes p_bt=t\otimes t\end{equation*}
Thus $\G\rtimes C_p$ consists of grouplikes, hence $H(\epsilon\otimes\epsilon,\tl)=\k(\G\rtimes C_p)$. In general, that is if $\ov{\chi_{\eta}}\ne\ov{\epsilon}$, $t$ can be twisted into a grouplike.

In the foregoing notation, let $f=\eta(t)$. We claim the element $ft$ is a grouplike in $H(\delta\eta,\tl)$. First, since $\eta:G\to\bu{(\D{C_p})},\,\eta(t)\in\D{G}$, and $ft$ makes sense. By \eqref{DeltaH}
\begin{equation*}\Delta_H(t)=\sum_{a,b}(\delta f)(a,b)p_at\o p_bt=\sum_{a,b}f(a)f(b)f(ab)^{-1}p_at\o p_bt\end{equation*}
Next apply $\Delta_H$ to the standard expansion $f=\sum_a f(a)p_a$. We get
\begin{equation*}\Delta_H(f)=\sum_a f(a)(\sum_{bc=a}p_b\o p_c)=\sum_{b,c}f(bc)p_b\o p_c\end{equation*}
All in all we have
\begin{align*}&\Delta_H(ft)=\Delta_H(f)\Delta_H(t)\\&=(\sum_{a,b}f(ab)p_a\o p_b)(\sum_{a,b}f(a)f(b)f(ab)^{-1}p_at\o p_bt)\\
&=\sum_{a,b}f(a)f(b)p_at\o p_bt=(\sum_af(a)p_a)t\o(\sum_{b}f(b)p_b)t=ft\o ft,\end{align*}
as needed.

Set $x=ft,\chi=\chi_{\eta}$ and observe that $x^p=\phi_p\bullet f=\chi$. We see $x$ is a unit in $H(\delta\eta,\tl)$. The action of $x$ on $\G$ by conjugation coincides with the action of $t$. Let $G(\ov{\chi},\bullet)$ be the subgroup of $H(\delta\eta,\tl)$ generated by $\G$ and $x$. Clearly $G(\ov{\chi},\bullet)/{\G}=C_p$, hence $G(\ov{\chi},\bullet)$ is an extension of $C_p$ by $\G$ associated to the datum $\{\ov{\chi},\bullet\}$. There $\ov{\chi}$ represents the cohomology class of $G(\ov{\chi},\bullet)$ as an element of $\mathrm{Opext}(C_p,\G,\bullet)$, since $H^2(C_p,\G,\bullet)=\G^{C_p}/N(\G)$.
From $|G(\ov{\chi},\bullet)|=\dim\,H(\delta\eta,\tl)$ we conclude $H(\delta\eta,\tl)=\k G(\ov{\chi},\bullet)$.\qed

It becomes apparent that we have reduced the isomorphism problem for Hopf algebras to same problem for groups $G(\ov{\chi},\bullet)$. We need to translate condition \eqref{isomorphismcocommutative} into a condition for the data $\{\ov{\chi},\bullet\}$ and $\{\ov{\omega},\circ\}$. In keeping with our convention we treat a coboundary $\delta\eta$ as an element of either $H^2_{\c}(\tl)$ or $H^2_{\c}/B^2_c(\tl)$.
\begin{Lem}\label{translation} Let $\phi:(\G,\bullet)\to(\G,\circ)$ be a $C_p$-isomorphism, and $\delta\eta\in H^2_{\c}(\tl),\delta\zeta\in H^2_c(\tl')$. Put $\ov{\chi}=\Phi\Theta(\delta\eta)$ and $\ov{\omega}=\Phi\Theta(\delta\zeta)$. Then $\delta\zeta=(\delta\eta).(\phi^*)^{-1}\;\text{iff}\;\phi(\ov{\chi})=\ov{\omega}$.\end{Lem}
\pf In the above statement we used the same letter $\phi$ for the induced isomorphism $\G^{C_p}/N(\G,\bullet)\to \G^{C_p}/N(\G,\circ)$. By definition of $(\phi^*)^{-1}$, $(\delta\eta).(\phi^*)^{-1}(a,b)=\eta(a\phi^*)\eta(b\phi^*)\eta((ab)\phi^*)^{-1}$. As $\eta(a\phi^*)=\phi(\eta)(a)$ by \eqref{dualphi}, it follows that $(\delta\eta).(\phi^*)^{-1}=\delta(\phi(\eta))$, hence the condition of the Lemma says $\delta\zeta=\delta(\phi(\eta))$. Applying $\Phi\Theta$ to the last equation, and using $C_p$-linearity of $\phi$ we derive $\phi_p\circ\zeta(t)N(\G,\circ)=\phi(\phi_p\bullet\eta(t))N(\G,\bullet)$, that is $\ov{\omega}=\phi(\ov{\chi})$. Since all steps of the proof are reversible, the proof is complete. \qed 

The final step of the proof is 

\begin{Prop}\label{maincocom}$G(\ov{\chi},\bullet)\simeq G(\ov{\omega},\circ)$ if and only if $\exists\, C_p$-isomorphism\\ $\phi:(\G,\bullet)\to(\G,\circ)$ such that $\phi(\ov{\chi})=\ov{\omega}.$\end{Prop}
\pf  The proof of the proposition will be carried out in steps. 

(1) We assume action `$\bullet$' to be nontrivial which is equivalent to assuming $H(\delta\eta,\tl)$ is a noncommutative algebra. We simplify notations by replacing $\G$ with $G$, and $\ov{\chi}$ by $\ov{a}$ where $a\in G^{C_p}$ and $\ov{a}=aN(G)$. An extension of $C_p$ by $G$ defined by some $\ov{a}$ and `$\bullet$' will be denoted by $G(\ov{a},\bullet)$. Recapping Lemma \ref{cocommutativegroup} we note that the group $G(\ov{a},\bullet)$ is generated by $G$ and an element $x\notin G$ such that $x^p=a$ and the action of $x$ in $G$ by conjugation coincides with the action of a generator of $C_p$. We note that if $\ov{a}=\ov{1}$ then $x$ can be chosen so that $x^p=1$. For, from $x^p=\phi_p\bullet b$ we have $(b^{-1}x)^p =1$. It follows each $G(\ov{1},\bullet)=G\rtimes C_p$. 

(2) The groups $G(\ov{1},\bullet)$ and $G(\ov{a},\circ)$ are nonisomorphic for any choice of $\ov{a}\ne\ov{1}$. There `$\circ$' denotes the action of $C_p$ in $G$ for the second group.

Suppose $\psi:G(\ov{1},\bullet)\to G(\ov{a},\circ)$ is an isomorphism. Let $x,y$ be elements of $G(\ov{1},\bullet)$ and $G(\ov{a},\circ)$,respectively with $x^p=1$ and $y^p=a$.  Were $\psi(g)=hy^k$ for some $g,h\in G$, we would have $1=\psi(g^p)=(\phi_p(y^k)\circ h)a^k=(\phi_p(y)\circ h)a^k$ contradicting to $\ov{a}\ne\ov{1}$. Thus $\psi(G)=G$, hence $\psi(x)=cy^k,\,c\in G,\,1\le k\le p-1$. But then the preceeding argument gives $1=\psi(x^p)=(\phi_p(y)\circ c)a^k$ whence $\ov{a}=\ov{1}$, a contradiction. Thus such $\psi$ does not exist.

(3) Here we show that $C_p$-modules $(G,\bullet^{\alpha})$ and $(G,\bullet)$ are isomorphic for every $\alpha\in A_p$. Let us write $G$ additively. Set $R=\mathbb{Z}_pC_p$ and $R_l=R/J_l$ where $J_l=\gen{(t-1)^l},1\le l\le p$ with $\gen{u}$ denoting the submodule generated by $u$. 
The action of $C_p$ in $G$ induces an $R$-module structure in $G$. Every indecomposable $R$-module is isomorphic to some $R_l$ with the action of $R$ by the left multiplication. Therefore the Krull-Schmidt decomposition of $(G,\bullet)$  
\begin{equation}\label{KSdecomposition}(G,\bullet)=B_1\oplus\cdots\oplus B_p,\end{equation}
consists of blocks $B_l=R_l^{m_l}$ of direct sums of modules $R_l$. There the action `$\bullet$' is taken to be the left multiplication. The sequence $\{m_l\}$ determines the isomorphism type of $(G,\bullet)$. An automorphism $\alpha:x\mapsto x^k,x\in C_p$ induces the automorphism of $R$ which sends $J_l$ to itself. Therefore $\alpha$ induces an automorphism of $R_l$, hence the $R$-module isomorphism $(R_l,\bullet)\simeq (R_l,\bullet^{\alpha})$.
This proves our claim.

(4) Here we prove the proposition  for groups $G(\ov{a},\bullet)$ and $G(\ov{b},\circ)$ with $\ov{a}\ne\ov{1}$ and $\ov{b}\ne\ov{1}$. We need only to show the necessary part, proof of sufficiency is fairly straghtforward.

Suppose $\psi: G(\ov{a},\bullet)\to G(\ov{b},\circ)$ is an isomorphism. Let $x,y$ be elements of those groups such that $x^p=a$ and $y^p=b$. By the argument used in (2) there holds: $\psi(G)=G$ and $\psi(x)=cy^k$ for some $c\in G$ and $1\le k\le p-1$. We derive the equality
\begin{equation}\label{twoactions} \psi(x\bullet g)=\psi(xgx^{-1})=y^kgy^{-k}=y^k\circ g.\end{equation}
Note \eqref{twoactions} shows the restriction $\phi=\psi|_G$ to be a $C_p$-isomorphism $\phi:(G,\bullet)\to (G,\circ^{\alpha})$ where $\alpha:x\mapsto x^k$. Furthermore $\psi(x)=cy^k$ implies $\phi(a)=\psi(x^p)=(\phi_p(y^k)\circ c)y^{pk}=(\phi_p\circ c)b^k$ as $\phi_p(y^k)=\phi_p(y)$. Say $\lambda:(G,\circ^{\alpha})\to (G,\circ)$ is a $C_p$-isomorphism guaranteed by part (3). One can see easily that $\phi_p(t)R=J_{p-1}$, hence $N(G,\circ)=J_{p-1}\circ G$ for any action $\circ$. Therefore $N(G,\circ^{\alpha})=J_{p-1}\circ^{\alpha}G=\alpha(J_{p-1})\circ G=J_{p-1}\circ G=N(G,\circ)$. As $\lambda$ is $R$-isomorphism, $\lambda(N(G,\circ^{\alpha}))=N(G,\circ)$, that is $\lambda(N(G,\circ))=N(G,\circ)$. Therefore as $b$ is a fixed point, so is $s=\lambda(b^k)$. Let $b_l,s_l$ be the $B_l$ components of $b,s$ from a decomposition \eqref{KSdecomposition} for $(G,\circ)$. Since $b$ and $s$ are fixed points, so are $b_l$ and $s_l$ for all $l$. Therefore they lie in the socle of $B_l$ and are simulteneously equal to $0$, or distinct from $0$, hence there is an automorphism of the socle mapping $s_l$ to $b_l$. Since each $B_l$ is a free $R_l$-module there exists a $C_p$-automorphism $\sigma_l$ such that $\sigma_l(s_l)=b_l$. It follows that there exists a $C_p$-automorphism $\sigma$ of $(G,\circ)$ with $\sigma(s)=b$.

Indeed, suppose $B_l=R_1^{(1)}\times\cdots\times R_l^{(m)}$. Let $\ov{1}=1+J_l$ be a generator of $R_l$ as a $C_p$-module, and $g_j$ be a copy of $\ov{1}$ in $R_l^{(j)}$. Since $b_l$ is a fixed point, $b_l\in (t-1)^{l-1}B_l$, hence $b_l=((t-1)^{l-1}(\sum k_jg_j),k_j\in\mathbb{Z}_p$. The mapping $\beta:g_1\mapsto k_1g_1\cdots k_mg_m,\,g_i\mapsto g_i,i>1$ extends to a $C_p$-automorphism with $\beta((t-1)^{l-1}g_1=b_l$. If $s_l$ is another element of $(t-1)^{l-1}B_l$, then then there exists a $C_p$-automorphism $\gamma:(t-1)^{l-1}g_1\to s_l$. But then $\beta\gamma^{-1}(s_l)=b_l$.

Since $\sigma$ commutes with the action of $C_p$, $\sigma\lambda(N(G,\circ))=N(G,\circ)$. The mapping $\tilde{\phi}=\sigma\lambda\phi$ is a $C_p$-automorphism $(G,\bullet)\to (G,\circ)$ with the property $\tilde{\phi}(a)=(\sigma\lambda(\phi_p\circ c))b$, hence $\tilde{\phi}(\ov{a})=\ov{b}$. This completes the proof of (4).

(5) We consider an isomorphism $\psi:G(\ov{1},\bullet)\to G(\ov{1},\circ)$. We need only to show the modules $(G,\bullet)$ and $(G,\circ)$ are isomorphic. 
Put $\mathbb{G}=G(\ov{1},\bullet)$ and $\mathbb{G}_1=G(\ov{1},\circ)$. For a group $F$ we let $\{\gamma_r(F)\}$ denote the lower central series of $F$ \cite{Ha}. A routine calculation yields

\noindent$\gamma_r(\mathbb{G})=(t-1)^{r-1}\bullet G$. One can see easily $\dim_{\mathbb{Z}_p}\displaystyle\gamma_l(\mathbb{G})/\gamma_{l+1}(\mathbb{G})=m_l+\cdots+ m_p$. It becomes evident that the multiplicities $m_j$ are determined by the lower central filtration. Since an isomorphism $\psi:\mathbb{G}\to\mathbb{G}_1$ induces isomorphism betweeen the lower central series in $\mathbb{G}$ and $\mathbb{G}_1$, the sequence $(m_1,\ldots,m_p)$ is an isomorphism invariant. This proves (5).

(6) It remains to settle the case of the trivial action. Now $G(\ov{a}):=G(\ov{a},\text{triv})$ is abelian. We have $N(G)=G^p=1$, hence $\ov{a}=a$. If $a=1, G(a)=G\times \gen{x}$ is an elementary $p$-group. Else, $a\ne 1,\,x^p=a$ which shows $x$ has order $p^2$. It is clear $G(1)\not\simeq G(a)$ for every $a\ne 1$. Furthermore, if $a\ne 1$, let $\ov{G}$ be a complement of $a$ in $G$. Evidently $G(a)=\ov{G}\times\gen{x}$, hence if $b\ne 1$ is another element of $G$, $G(a)\simeq G(b)$. On the other hand, for every $a,b\ne 1$ there is an automorphism $\phi$ of $G$ with $\phi(a)=\phi(b)$. This completes the proof of the proposition.\qed

For $p=2$ the isomorphism theorems yield
\begin{Cor}\label{elem2groups} Let $G$ be an elementary $2$-group. Then there is a bijection between the orbits of $\A$ in $H^2_c(\tl)$ and the isoclasses of extensions in $\E$.
\end{Cor}


\section{Commutative Extensions}\label{commutative}

An algebra $H(\tau,\tl)$ is commutative iff the action `$\tl$' is trivial. Below we omit the symbol `$\tl$' and write $H(\tau)$ for $H(\tau,\text{triv})$. Every commutative finite-dimensional Hopf algebra over an algebraically closed field is of the form $\D{L}$ (\cite{LR}, \cite[2.3.1]{Mo}) for some finite group $L$. We will identify groups $L$ appearing in that formula for algebras $H(\tau)$. It is convenient to introduce the group $\mathrm{Cext}(G,C_p)$ of central extensions of $C_p$ by $G$ \cite{BT}.
\begin{Prop}\label{commextensions} The group $\mathrm{Ext}_{[\text{triv}]}(\k C_p,\k^{G})$ of equivalence classes of commutative extensions is isomorphic  to the group $\mathrm{Cext}(G,C_p)$. The isomorphism is given by $H(\tau)\leftrightarrows \D{G(\tau)}$ where $G(\tau)$ is the central extension defined by the $2$-cocycle $\tau$.
\end{Prop} 
\pf In one direction, pick $\tau\in H^2_c(\text{triv})$. We have by Proposition \ref{cocrossedproduct1} for the trivial action that for every $a\in G$ 
\begin{equation*}\Delta_{H^*}(\ov{a})=\sum_{x,y\in C_p}\ov{a}p_x\o \ov{a}p_y=(\ov{a}\o\ov{a})(\sum_xp_x\o\sum_yp_y)=\ov{a}\o\ov{a}\end{equation*}
Thus $\ov{a}$ is a grouplike for every $a\in G$. Let $\theta$ be a generator of $\w{C}_p$. Since $\D{C_p}$ is a Hopf subalgebra of $H^*(\tau)$, $\theta$ is a grouplike of $H^*(\tau)$. We see that the set $G(\tau)=\{\ov{a}\theta^i|a\in G,0\le i\le p-1\}$ consists of grouplikes. Moreover, $|G(\tau)|=\dim\,H^*(\tau)$, hence $G(\tau)$ is a basis of $H^*(\tau)$. Therefore $G(\tau)=G(H^*(\tau))$ is a group and $H^*(\tau)=\k G(\tau)$, whence $H(\tau)=\D{G(\tau)}$.

Since $\tau$ is a Hopf $2$-cocycle it satisfies \eqref{hopfcocycle} which for the trivial action of $C_p$ turns into $\tau(xy)=\tau(x)\tau(y)$. Thus  $\tau(a,b):C_p\to\bu{\k}$ is a character for any choice of $a,b\in G$. We see that $\tau:G\times G\to\w{C_p}$ is a $2$-cocycle of $G$ in $\w{C_p}$, hence $G(\tau)$ is a central extension of $G$ by $\w{C_p}$ defined by $\tau$. It remains to notice that the subgroup $B^2_c$ consists of coboundaries in the group  $Z^2(G,\w{C}_p)$. For, by Definition \ref{Hopfcocycles} $\delta\eta\in B^2_c$ iff  $\eta$ satisfies the condition \eqref{hopf1cocycle}, hence $\eta: G\to \w{C}_p$.

The opposite direction is trivial.\qed

For calculations of orbits of $\g(\tl)$ in $H^2_c(\tl)$ we prefer to use a much smaller space 
\begin{equation}\label{classifyingspace}\mathbb{X}(\tl)=H^2(C_p,\G,\tl)\times H^2_N(G,\bu{\k}).\end{equation}
We make several remarks regarding $\mathbb{X}(\tl)$. Let $C_p$ act on $\k^G$ by $\bullet$ with the induced action $\tl$ on $G$. Then $C_p$ acts on $Z^2(G,(\k^{C_p})^{\bullet})$ via $x\bullet\tau(y,a,b)=\tau(y,a\tl x,b\tl x),x\in C_p,a,b\in G$. Recall $Z^2_N(G,\k^{\bullet})$ is the subgroup of $Z^2(G,\k^{\bullet})$ of all $s\in Z^2(G,\k^{\bullet})$ subject to $\phi_p\bullet s=1$, where $\phi_p\bullet s=\prod(t^i\bullet s)$, and  $\phi_p=1+t+\cdots+t^{p-1}$. The next lemma is a strengthening of Lemma \ref{actiononcocycles}.
\begin{Lem}\label{actiononX} {\rm(i)} For every $\lambda\in I(\tl,\tl')$ the mapping $s\mapsto s.\lambda$ is a $C_p$ isomorphism $\mathbb{X}(\tl)\simeq\mathbb{X}(\tl')$.

{\rm(ii)} $\mathbb{X}(\tl)$ is invariant under action by elements of $I(\tl,\tl^{\alpha})$ for every $\alpha\in A_p$.
\end{Lem}
\pf {\rm(i)} Let $\circ$ denote the action of $C_p$ on $Z^2(G,\k^{\bullet})$ induced by $\tl'$, i.e. $(x\circ s(a,b)=s(a\tl' x,b\tl' x)$. We need to show  $\phi_p\circ(s.\lambda)=1$. Since $\lambda^{-1}\in I(\tl',\tl)$ there holds $(a\tl' x).\lambda^{-1}=a\lambda^{-1}\tl x$. Therefore $t^i\circ(s.\lambda)=(t^i\bullet s).\lambda$ as the following calculation shows:
\begin{align}\label{intertwining} t^i\circ(s.\lambda)(a,b)&=(s.\lambda)(a\tl' t^i,b\tl' t^i)=s((a\tl' t^i).\lambda^{-1},(b\tl't^i).\lambda^{-1})\\
                                       &=s(a.\lambda^{-1}\tl t^i,b.\lambda^{-1}\tl t^i)=(t^i\bullet s).\lambda(a,b)\nonumber.\end{align}
We conclude that $\lambda$ induces a $C_p$-linear map $\x(\tl)\to\x(\tl')$, and also 
\begin{equation*}\phi_p\circ(s.\lambda)=\prod(t^i\circ(s.\lambda))=(\prod(t^i\bullet s)).\lambda=1\end{equation*}

{\rm(ii)} We must show the equality $\phi_p\bullet(s.\lambda)=1$ for every $s\in Z^2_N(\tl)$. By part (i) with $\circ=\bullet^{\alpha}$ there holds $\phi_p\bullet^{\alpha}(s.\lambda)=1$. Assuming $\alpha:x\mapsto x^k$ and noting $\phi_p(t^k)=\phi_p(t)$ this gives \begin{equation*}1=\phi_p\bullet^{\alpha}(s.\lambda)=\phi_p(t^k)\bullet(s.\lambda)=\phi_p\bullet(s.\lambda)\end{equation*} 
Since $\lambda$ sends $\ker\Phi$ in $Z^2_N(\tl)$ to $\ker\Phi$ in $Z^2_N(\tl^{\alpha})$ the proof is complete.\qed

We turn $\mathbb{X}(\tl)$ into a $\g(\tl)$-module by transfering its action from $H^2_c(\tl)$ to $Z^2_N(\tl)$ along $\Theta$. Pick $\omega_{\lambda,\alpha}$ and suppose $\alpha^{-1}:x\mapsto x^l,x\in C_p$. For $s\in Z^2_N(G,\bu{\k})$ we set
\begin{equation}\label{derivedaction}s.\omega_{\lambda,\alpha}=(\phi_l\bullet s).\lambda.\end{equation} 
Let $\Theta_*$ be the isomorphism of Corollary \ref{basicisic}.
\begin{Lem}\label{gisomorphism} {\rm(i)}. For every prime and any action `$\tl$' isomorphism $\Theta_*:H^2_c(\tl)\simeq Z^2_N(G,\bu{\k})/\ker\Phi$ is  $\g(\tl)$-linear.\\  
{\rm(ii)}. For every odd $p$ the isomorphism $Z^2_N(G,\bu{\k})/\ker\Phi\simeq H^2(C_p,\G)\times H^2_N(G,\bu{\k})$ is $\g(\tl)$-linear.\end{Lem}
\pf (i). For every $\tau\in Z^2_c(\tl)$ we have by \eqref{componentsoftau} 
\begin{align*}&\Theta(\tau.\omega_{\lambda,\alpha})=(\tau.\omega_{\lambda,\alpha})(t)=(\tau.\lambda)(t^l)=\phi_l\bullet^{\alpha}((\tau.\lambda)(t))=\phi_l\bullet^{\alpha}(\tau(t).\lambda)\\&(\text {by}\, \eqref{intertwining})=(\phi_l\bullet\tau(t)).\lambda=\Theta(\tau).\omega_{\lambda,\alpha}\end{align*}
This equation demonstrates that \eqref{derivedaction} turnes $Z^2_N(\tl)$ into a $\g(\tl)$-module. It is immediate that $B^2_c(\tl)$ is a $\g(\tl)$-subgroup of $Z^2_c(\tl)$. By Lemma \ref{structureofcoboundaries}(ii) $\ker\Phi$ is a $\g(\tl)$-subgroup, which proves part (i).

(ii). For an odd $p$ splitting \eqref{Z^2adm} is carried out by the mapping $s\mapsto s\un{a}(s^{-2})\times\un{a}(s^2)$ which is clearly a $\g(\tl)$-map. It remains to note that homomorphism $\Phi$ is also a $\g(\tl)$-map.\qed

The next result gives the number of isotypes of commutative extensions in $\mathrm{Ext}(\k C_p,\D{G})$ for odd primes $p$ and elementary $p$-groups. By Proposition \ref{commextensions} this is also the number of nonisomorphic groups in $\mathrm{Cext}(G,C_p)$ which can be derived by group-theoretic methods. Our proof avoids group theory and sets up framework for generalization to non-elementary abelian groups. However, it does not cover $2$-groups.
\begin{Prop}\label{commisotypes} Let $G$ be a finite elementary $p$-group of order $p^n$ for an odd $p$. There are $\lfloor\frac{3n+2}{2}\rfloor$ isotypes of commutative Hopf algebra extensions of $\k C_p$ by $\D{G}$ for any odd prime $p$.\end{Prop}

\pf Write $\g$ for $\g(\text{triv})$ and likewise $\At=\mathbb{A}(\text{triv})$. By Theorems \ref{mainthm}, \ref{cocommutativeclasses} and the preceeding lemma  the isotypes in question are in a bijection with the set of $\g$-orbits in $\x(\t)$. We observe that for every trivial $C_p$-module $M$, $N(M)=M^p$. Therefore, if $M$ has exponent $p$ $N(M)=1$ and $M_N=M$. From these remarks, in view of $\G$ and $\mathrm{Alt}(G)$ having exponent $p$, we derive by Lemma \ref{gisomorphism} a $\g$-isomorphism
\begin{equation}\label{commcohomology} H^2_c(\text{triv})\simeq \G\times\text{Alt}(G).\end{equation}
By definition $\g=\At\times A_p$ with elements $\alpha_k:x\mapsto x^k$ of $A_p$ acting on $\x(\t)$ by \eqref{derivedaction} as 
\begin{equation*} z.\alpha_k=\phi_k\bullet z=z^k\;\text{as\,$\bullet$\,is trivial}\end{equation*}
We proceed to description of orbits of $\At$ in $H^2_c(\text{triv})$. This description will show that every $\At$-orbit is closed to action by $A_p$, hence $\g$- and $\At$-orbits coincide. We switch to the additive notation in our treatment of $\G\times\text{Alt}(G)$ and view the latter as a vector space over the prime field $\mathbb{Z}_p$. We note $\mathbb{A}$ acts in $\G\oplus\text{Alt}(G)$ componentwise. For a $(\chi,\beta)\in\G\times\text{Alt}(G)$ we let $(\chi,\beta)\At$ denote the orbit of $(\chi,\beta)$. For every $\chi\in\G$ we let $K_{\chi}$ denote $\ker\chi$.

Classification of orbits of $\mathbb{A}$ in $\G\oplus\text{Alt}(G)$ relies on the theory of symplectic spaces. A symplectic space $(V,\beta)$ is a vector space with an alternating form $\beta$. We need a structure theorem for such spaces (see e.g. \cite{La}). For a subspace $X$ of $V$ we denote by $X^{\perp}$ the subspace consisting of all $v\in V$ such that $\beta(x,v)=0$ for all $x\in X$. We call $V^{\perp}$ the radical of $\beta$ and denote it by $\text{rad}\,\beta$. We say that two elements $x,y\in V$ are orthogonal if $\beta(x,y)=0$. We call subspaces $X,Y$ orthogonal if $X\subset Y^{\perp}$, and we write $X\perp Y$ in this case. A symplectic space $(V,\beta)$ is an orthogonal sum of subspaces $U_1,\ldots, U_k$ if $V=U_1\oplus\cdots\oplus U_k$ and $U_i\perp U_j$ for all $i\ne j$. We use the symbol $V=U_1\perp\cdots\perp U_k$ to denote an orthogonal decomposition of $(V,\beta)$. A hyperbolic plane $P$ is a $2$-dimensional subspace of $V$ such that $\beta(x,y)=1$ relative to a basis $\{x,y\}$ of $P$. The elements $x,y$ are called a hyperbolic pair. The structure theorem states 
\begin{equation}\label{symplectic} V=P_1\perp\cdots\perp P_r\perp\text{rad}\,\beta .\end{equation}
We call this splitting of $V$ a complete orthogonal decomposition of $(V,\beta)$. The number $r$ will be refered to as the width of $\beta$, denoted by $w(\beta)$. 

Let us agree to write $\text{Alt}_r(G)$ for the set of alternating forms of width $r$. We identify $\G\oplus\text{Alt}_0(G)$ with $\G$. The set $\G\oplus\text{Alt}_r(G)$ is visibly stable under $\mathbb{A}$. Using \eqref{symplectic} one can see easily $\mathbb{A}$ acts transitively on $\text{Alt}_r(G)$ for every $r$. It folows that every orbit of $\mathbb{A}$ lies in $\G\oplus\text{Alt}_r(G)$ for some $r$. A refinement of $\G\oplus\text{Alt}_r(G)$ gives all orbits of $\mathbb{A}$.
\begin{Prop}\label{orbits} The orbits of $\mathbb{A}$ in $\G\oplus\text{Alt}_r(G)$ are as follows:

{\rm (i)} $\{0\}$ and $\G\setminus\{0\}$ if $r=0$;

{\rm(ii)} $\{(0,\beta)\},\;\{(\chi,\beta)|\text{rad}\,\beta\subset K_{\chi}\}$,\;  $\{(\chi,\beta)|\text{rad}\,\beta\not\subset K_{\chi}\}$ for every $1\le r\le\lfloor n/2\rfloor$ for an odd $n$, or $1\le r< n/2$ for an even $n$, where $\beta$ runs over $\mathrm{Alt}_r(G)$;

{\rm(iii)} $\{(0,\beta)\}$ and $\{(\chi,\beta)|0\ne \chi\}$, where $\beta$ runs over $\mathrm{Alt}_{n/2}$ if $n$ is even.
\end{Prop}
\pf Pick $\chi\in\G$. Let $x\in G$ be an element, unique modulo $K_{\chi}$, such that $\chi(x)=1$ where $1$ is the unity of $\mathbb{Z}_p$. Clearly $\chi$ is uniquely determined by a pair $(K_{\chi},x)$. For any $\lambda\in\G$ with an associated pair $(K_{\lambda},z)$ the equality $\chi.\phi=\lambda,\phi\in\mathbb{A}$ holds iff $K_{\chi}\phi=K_{\lambda}$ and $x\phi=k+z$ for some $k\in K_{\lambda}$. On the other hand, $\beta,\gamma\in\mathrm{Alt}(G)$ are related by $\beta.\phi=\gamma$ iff there is a decomposition \eqref{symplectic} of $(V,\beta)$ satisfying $P_i\phi$ is a hyperbolic plane for $\gamma$ for all $i$ and $(\text{rad}\,\beta)\phi=\text{rad}\,\gamma$.

Let $\gen{X}$ denote the subspace spanned by a subset $X\subset G$. Note $\mathbb{A}$ acts transitively on the set of pairs $(L,x)$  such that $G=L\oplus\gen{x}$ whence $\mathbb{A}$ acts transitively on $\G\setminus\{0\}$ which proves (i). (iii) is a special case of the second set in (ii) as $\text{rad}\,\beta=0$ for every $\beta$ of width $n/2$.

We take up part (ii). First, we show that all sets there are $\mathbb{A}$ invariant. It suffices to consider the property $\text{rad}\,\beta\subset\ker\chi$ of $(\chi,\beta)$. Let $(\lambda,\gamma)=(\chi,\beta).\phi$. Then $\beta.\phi=\gamma$ implies $(\text{rad}\,\beta)\phi=\text{rad}\,\gamma$. The second condition $\lambda=\chi.\phi$ yields the equality $K_{\chi}\phi=K_{\lambda}$. Therefore $\text{rad}\,\gamma=(\text{rad}\,\beta)\phi\subset K_{\chi}\phi=K_{\lambda}$.

Second, we prove that each set in (ii) is a single orbit. We begin with $\{(\chi,\beta)|\text{rad}\,\beta\subset K_{\chi}\}$. Let $\beta_{\chi}$ denote the restriction of $\beta$ to $K_{\chi}$. Take a complete orthogonal decomposition $(K_{\chi},\beta)=P_1\perp\cdots\perp P_m\perp\text{rad}\,\beta_{\chi}$. Suppose $w(\beta)=r$. Then $\dim\text{rad}\,\beta=n-2r$ while $\dim\text{rad}\,\beta_{\chi}=n-1-2m$. Since $\text{rad}\,\beta\subset\text{rad}\,\beta_{\chi}$ it follows that $m< r$. Therefore $\dim\text{rad}\,\beta_{\chi}\ge n-1-2(r-1)=n-2r+1>\dim\text{rad}\,\beta$. Further select an $x\in\text{rad}\,\beta_{\chi}\setminus\text{rad}\,\beta$ and some $y\notin K_{\chi}$. Notice $x$ is not orthogonal to $y$ for else, as $G=K_{\chi}\oplus\gen{y}$, $x\in\text{rad}\,\beta$, a contradiction. Let $R=\text{rad}\,\beta_{\chi}$. Since $\beta(x,y)\ne0$ the functional $\beta(-,y):R\to\mathbb{Z}_p,r\mapsto\beta(r,y),r\in R$ is nonzero. Therefore $R$ splits up as $R=\gen{x}\oplus\ker\beta(-,y)$. It is clear that $\ker\beta(-,y)=\text{rad}\,\beta$ which implies $\dim R=n-1-2m=n-2r+1$ hence $m=r-1$. Put $P=P_1\oplus\cdots\oplus P_{r-1}$ and observe that the restriction of $\beta$ to $P$ is nondegenerate forcing $P\cap P^{\perp}=0$. From this we obtain a decomposition $G=P\perp P^{\perp}$. Since the number $w(\beta)$ is an invariant of decompositions \eqref{symplectic} and $\text{rad}\,\beta|_{P^{\perp}}=\text{rad}\,\beta$, we conclude that there is a single hyperbolic plane $P_r$ such that $P^{\perp}=P_r\perp\text{rad}\,\beta$. In consequence $R\cap P^{\perp}=(R\cap P_r)\perp\text{rad}\,\beta$ as $\text{rad}\,\beta\subset R$. Let $u,v$ be a hyperbolic pair in $P_r$ with $u\in R$. Then $v\notin K_{\chi}$, for otherwise $K_{\chi}\supset P_r$, and then $K_{\chi}=G$. We see that $G$ decomposes in two ways
\begin{align}&(G,\beta)=K_{\chi}\perp\gen{v},\,\text{and}\label{orthodecomposition1}\\
&(G,\beta)=P_1\perp\cdots\perp P_{r-1}\perp P_r\perp\text{rad}\,\beta\label{orthodecomposition2}\end{align}
with $K_{\chi}=P_1\oplus\cdots P_{r-1}\oplus\gen{u,\text{rad}\,\beta}$ and $P_r=\gen{u,v}$. Pick another element $(\lambda,\gamma)$ of the set. By \eqref{orthodecomposition1},\eqref{orthodecomposition2} $(G,\gamma)=Q_1\perp\cdots\perp Q_r\perp\text{rad}\,\gamma=K_{\lambda}\perp\gen{w}$ with $Q_r=\gen{z,w}$ and $K_{\lambda}=Q_1\oplus\cdots\oplus Q_{r-1}\oplus\gen{z,\text{rad}\,\gamma}$. Then any automorphism $\phi$ sending $P_i$ to $Q_i$ for $i=1,\ldots,r-1$, $\text{rad}\,\beta$ to $\text{rad}\,\gamma$, and $u\mapsto z,v\mapsto w$ carries $(\chi,\beta)$ to $(\lambda,\gamma)$.

Finally we show that any set  $\{(\chi,\beta)|\text{rad}\,\beta\not\subset K_{\chi}\}$ is a single orbit. Pick $(\chi,\beta)$ from the set, and let $(K_{\chi},\beta)=P_1\perp\cdots\perp P_m\perp\text{rad}\,\beta_{\chi}$ be a complete orthogonal decomposition of $(K_{\chi},\beta)$. Select some $y\in\text{rad}\,\beta\setminus K_{\chi}$. Then $G=K_{\chi}\perp\gen{y}$, hence $\text{rad}\,\beta=\text{rad}\,\beta_{\chi}\oplus\gen{y}$. We see $(G,\beta)=P_1\perp\cdots\perp P_m\perp (\text{rad}\,\beta_{\chi}\oplus\gen{y})$ is a complete orthogonal decomposition of $G$. It becomes evident that $m=r$ and for any other pair $(\lambda,\gamma)$ from the set $(\chi,\beta).\phi=(\lambda,\gamma)$ for some $\phi\in\mathbb{A}$. A simple count of the number of orbits yields the formula.\qed 

\section{Some Extensions of Dimension $p^4$}\label{p^4}

In this section $p$ is an odd prime and $G$ is an abelian group of order $\le p^3$.  When $|G|=p$, any action of $C_p$ on itself is trivial, hence $\mathrm{Ext}(\k^{C_p},\k C_p)=\mathrm{Ext}_{[\t]}(\k^{C_p},\k C_p)$ which by Theorem \ref{commisotypes} has two isoclasses. Moreover, as $\text{Alt}(C_p)=1$ every $H\in\mathrm{Ext}_{[\t]}(\k C_p,\k^{C_p})$ is cocommutative. By Proposition \ref{maincocom} part (6) $H=\k C_{p^2}$ or $H=\k(C_p\times C_p)$, and we derive part of \cite[Theorem 2]{Mas2}. The case $|G|=p^2$, also due to A. Masuoka \cite{Mas1}, will be dealt with below as  specialization of a more general theory for $|G|=p^3$. 

We assume that, unless stated otherwise, $|G|=p^3$. In the additive notation $G=\Z_p^3$ or $G=\Z_{p^2}\oplus\Z_p$, and the theory splits into two parts.

(A) Suppose $G=\Z_p^3$. There are up to isomorphism two nontrivial $\Z_pC_p$-module structures on $G$. Let $R_i=\Z_pC_p/\gen{(t-1)^i},0\le i\le p-1$. Then either $G\simeq R_2\oplus R_1$ or $G\simeq R_3$. Before proceeding to cases we make a notational change. We write $\alpha_k$ for the mapping $x\mapsto x^k,x\in C_p$ and $\tl^k$ for $\tl^{\alpha_k}$.

(I) Suppose $G\simeq R_2\oplus R_1$, and let $\tl$ be the action of $C_p$ on $G$ composed of regular actions of $C_p$ on $R_2$ and $R_1$. We aim to prove
\begin{Thm}\label{isoclassesAI} $\mathrm{Ext}_{[\tl]}(\k^{C_p^3},\k C_p)$ contains $2p+11$ isoclasses of extensions.\end{Thm}
We break up proof in steps.
 
 (1) Here we compute $\X$. Select a basis $\{e,g,f\}$ for $G$ where $\{e,f\}$ span $R_2$, and $R_1=\Z_pg$ with the action
\begin{equation}\label{taction}e\tl t=e+f, g\tl t=g, f\tl t=f.\end{equation}
Clearly the matrix $T$ of $t$ in that basis is $T=\begin{pmatrix}1&0&1\\0&1&0\\0&0&1\end{pmatrix}$. Let $\{e^*,g^*,f^*\}$ be the dual basis for $\G$, and   let $\wedge$ denote the multiplication in the Grassman algebra over $\G$.  We fix a basis $\{e^*\we g^*,e^*\we f^*,g^*\we f^*\}$ for $\G\we\G$, hence a basis $\{e^*,g^*,f^*,e^*\we g^*,e^*\we f^*,g^*\we f^*\}$ for $\G\oplus \G\we\G$. We refer to the above bases as standard. 
\begin{Prop}\label{structureX1} $\X=\G^{C_p}\oplus\G\wedge\G=\gen{e^*,g^*}\oplus\G\wedge\G$.\end{Prop} 
\pf Recall $\X= \G^{C_p}/N(\G)\oplus \text{Alt}_N(G)$. We use the well known identification  $\text{Alt}(G)=\G\we\G$. One can see easily that action of $t$ in $\G$ is described by $T^{\Tr}$ in the standard basis of $\G$. By general principles \cite[III,8.5]{Bou} the matrix of $t$ in the standard basis of $\G\we\G$ is $T^{\Tr}\we T^{\Tr}=\begin{pmatrix}1&0&0\\0&1&0\\-1&0&1\end{pmatrix}$.
It follows that $(t-1)^{p-1}\bullet\G=0$ and $(t-1)^{p-1}\bullet\G\we\G=0$, that is $N(\G)=0$ and $(\G\we\G)_N=\G\we\G$. Further, one can see easily $\G^{C_p}=\gen{e^*,g^*}$.\qed

(2) Group $\A$. Identifying $\phi\in\A$ with its matrix $\F$ one has $\F\in\A$ iff $\F T=T\F$. This condition leads up to a determination of $\A$, viz.
\begin{equation}\label{A1}\A=\left\{\F|\F=\begin{pmatrix}a_{11}&a_{12}&a_{13}\\0&a_{22}&a_{23}\\0&0&a_{11}\end{pmatrix},\,a_{ij},\in\Z_p,\,a_{11}a_{22}\ne 0\right\}\end{equation}

(3) Orbits of $\A$ in $\X$. Since $\A$ acts on $\G$ by $(v^*.\phi)(a)=v^*(a.\phi^{-1}),v^*\in\G$, the matrix of $\phi$ in the standard basis for $\G$ is $(\F^{-1})^{\Tr}$. We prefer to use coordinates $u,v,q,r,s$ for $\F^{-1}$ where $u=a_{11}^{-1},v=a_{22}^{-1}$ and $q=ua_{12},r=ua_{13}, s=va_{23}$. A routine calculation gives
\begin{equation}\label{Finverse}(\F^{-1})^{\Tr}=\begin{pmatrix}u&0&0\\-vq&v&0\\u(qs-r)&-us&u\end{pmatrix}\end{equation}
We treat the tuple $(u,v,q,r,s)$ as coordinates of either $\phi$ or $\F$. On general principles \cite[III,8.5]{Bou} the matrices of $\phi$ in the standard bases for $\G$ and $\G\we\G$ are $(\F^{-1})^{\Tr}$ and $(\F^{-1})^{\Tr}\we(\F^{-1})^{\Tr}$, respectively. For $\F$ defined by $(u,v,q,r,s)$ the result is 
\begin{equation}\label{FwedgeF}(\F^{-1})^{\Tr}\we(\F^{-1})^{\Tr}=\begin{pmatrix}uv&0&0\\-u^2s&u^2&0\\uvr&-uvq&uv\end{pmatrix}.\end{equation}
Regarding $\Z_p$ as field we let $\zeta$ denote a generator of $\bu{\Z_p}$. We observe a simple lemma
\begin{Lem}\label{partialorbits} There are two and four nonzero orbits of $\A$ in $\G^{C_p}$ and $\G\we\G$, respectively.\end{Lem}
\pf We assign a vector $(a_1,a_2)$ to the element $a_1e^*+a_2g^*$ of $\G^{C_p}$ and likewise $(b_1,b_2,b_3)$ to $b_1e^*\we g^*+b_2e^*\we f^*+b_3g^*\we f^*$. By \eqref{Finverse} $(0,1).\A=\{(a_1,a_2|a_2\ne 0\}$, and $(1,0).\A=\{(a_1,0)|a_1\ne 0\}$. This proves the first claim.

Similarly, using \eqref{FwedgeF} one can derive readily the equalities
\begin{equation*}(0,0,1).\A=\{(b_1,b_2,b_3)|b_3\ne 0\},\,(1,0,0).\A=\{(b_1,0,0)|b_1\ne 0\}.\end{equation*}
However, the set $\{(b_1,b_2,0)|b_2\ne 0\}$ is union of two orbits. Namely, if $b_2\in{\bu{\Z_p}}^2$, then by \eqref{FwedgeF} $(b_1,b_2,0)\in(0,1,0).\A$. But if $b_2\notin{\bu{\Z_p}}^2$, then $b_2\in\zeta\Z_p^2$ which implies $(b_1,b_2,0)\in(0,\zeta,0).\A$.\qed

We introduce notation
\begin{equation*}\Omega'_0=\{(0,0)\}, \Omega'_1=(1,0).\A,\Omega'_2=(0,1).\A,\;\text{and}\,\end{equation*} \begin{align*}\Omega''_0&=\{(0,0)\},\Omega''_1=(1,0,0).\A,\Omega''_{2,0}=(0,1,0).\A,\\\Omega''_{2,1}&=(0,\zeta,0).\A,\Omega''_3=(0,0,1).\A.\end{align*} Some of products $\Omega'_i\times\Omega''_j$ are orbits itself. We list those that are in
\begin{Lem}\label{splitorbits} The following sets are orbits
\bi
\item[(0)] $\Omega'_0\times\Omega''_j$ and $\Omega'_i\times\Omega''_0,\;j=0,1,(2,0),(2,1),3,\,i=1,2$.
\item[(1)] $\Omega'_1\times\Omega''_1$ and $\Omega'_1\times\Omega''_3$.
\item[(2)] $\Omega'_2\times\Omega''_1,\,\Omega'_2\times\Omega''_{2,0}$ and $\Omega'_2\times\Omega''_{2,1}$
\ei\end{Lem}
\pf Vectors $(a_1,a_2)$ and $(b_1,b_2,b_3)$ give rise to a concatenated vector $(a_1,a_2;b_1,b_2,b_3)$. The claim is that concatenating generators of $\Omega'_i,\Omega''_j$ for $i,j$ as in the Lemma we get a generator for $\Omega'_i\times\Omega''_j$. We give details for $\Omega'_1\times\Omega''_3$. Other cases are treated similarly. Combining \eqref{Finverse} with \eqref{FwedgeF} we obtain 
\begin{equation*}(1,0;0,0,1).\A=\{(u,0;uvr,-uvq,uv)|uv\ne 0,q,r,s\,\text{arbitrary}\}\end{equation*}
Now for every element $(a_1,0;b_1,b_2,b_3)\in\Omega'_1\times\Omega''_3$ the equations $u=a_1, uv=b_3,uvr=b_1,-uvq=b_2$ are obviously solvable, which completes the proof.\qed

We pick up $p-1$ additional orbits in 
\begin{Lem}\label{extraorbits} Each set $\Omega'_1\times\Omega''_{2,i},i=0,1$ is union of $(p-1)/2$ orbits.\end{Lem}
\pf Say $i=0$. By definition $\Omega'_1\times\Omega''_{2,0}=\{(a_1,0;b_1,b_2,0)|a_1\in\bu{\Z_p},b_2\in{\bu{\Z_p}}^2,b_1\,\text{arbitrary}\}$. For every $m\in{\bu{\Z_p}}^2$ we let\\$z_m=(1,0;0,m,0)$. By \eqref{Finverse} and \eqref{FwedgeF} $z_m.\phi=(u,0;-u^2sm,u^2m,0)$ where $u,s$ are among parameters of $\phi$. It is immediate that $|z_m\A|=(p-1)p$, and one can verify equally directly that $z_m.\A\cap z_n.\A=\emptyset$ for $m\ne n$. Since   $|\Omega'_1\times\Omega''_{2,0}|=\frac{(p-1)}{2}(p-1)p$ this case is done. For $i=1$ one should take $z'_m=(1,0;0,\zeta m,0)$.\qed

We summarize
\begin{Lem}\label{Aorbits}  There are $2p+11$ orbits of $\A$ in $\X$.\end{Lem}
\pf The previous two lemmas give $p+11$ orbits. The rest will come from splitting of the remaining set $\Omega'_2\times\Omega''_3$. The latter is defined as $\{(a_1,a_2;b_1,b_2,b_3)|a_2,b_3\in\bu{\Z_p},a_1,b_1,b_2\,\text{arbitrary}\}$. For every $k\in\Z_p$ we define $w_k=(k,1;0,0,1)$. Again by \eqref{Finverse} and \eqref{FwedgeF} we have
\begin{equation*}w_k.\A=\{(uk-vq,v;uvr,-uvq,uv)\}.\end{equation*}
There $u,v$ run over $\bu{\Z_p}$ and $r,q$ run over $\Z_p$. One can see easily that $|w_k.\A|=(p-1)^2p^2$. Furthermore, we claim that $w_k.\A\cap w_l.\A\\=\emptyset$ for $k\ne l$.

For, suppose 
\begin{equation*}(uk-vq,v;uvr,-uvq,uv)=(u'l-v'q',v';u'v'r',-u'v'q',u'v')\end{equation*}
for some $(u,v,q,r)\,\text{and}\, (u',v',q',r')$. Then $v=v'$ and $uv=u'v'$ give $u=u'$. This implies $q=q',r=r'$, and finally $uk=ul$ yields $k=l$, a contradiction. We conclude that $|\cup_{0\le k\le p-1}w_k.\A|=p^3(p-1)^2$. As this is the number of elements in $\Omega'_2\times\Omega''_3$, the proof is complete.\qed 

(4) Orbits of $\g(\tl)$. By definition $\g(\tl)$ is generated by $\A$ and a select set $\{\lambda_k|2\le k\le p-1\}$ of automorphisms of $G$. There $\lambda_k\in I(\tl,\tl^k)$, and by \eqref{modulemap} $\lambda\in I(\tl,\tl^k)$ iff its matrix $\Lambda$ satisfies
\begin{equation}\label{intertwiner}T\Lambda=\Lambda T^k.\end{equation}
Set $\Lambda_k=\text{diag}(1,1,k)$ (that is the diagonal matrix with entries $1,1,k$), and observe that $\Lambda_k$ satisfies \eqref{intertwiner}. We denote by $\lambda_k$ the automorphism whose matrix is $\Lambda_k$, and we set $\omega_k=\lambda_k\alpha_k^{-1}$. We move on to calculation of matrices of automorphisms $\omega_k$. We set $l=k^{-1} \pmod p$.  
\begin{Lem}\label{gaction} Action of $\omega_k$ is described by
\begin{align*} e^*.\omega_k&=le^*,\,g^*.\omega_k=lg^*\\
e^*\we g^*.\omega_k&=le^*\we g^*\\e^*\we f^*.\omega_k&=l^2e^*\we f^*\\g^*\we f^*.\omega_k&=-\binom l2 e^*\we g^*+l^2g^*\we f^*\end{align*}\end{Lem}
\pf By \eqref{derivedaction} for $\tau\in\X$, $\tau.\omega_k=(\phi_l\bullet\tau).\lambda_k$. For $\tau=e^*,g^*,e^*\we g^*,e^*\we f^*$ $\phi_l\bullet\tau=l\tau$ as such $\tau$ is fixed by $C_p$. Because $(t-1)^2\bullet\G\we\G=0$ we expand $\phi_l$ in powers of $t-1$, namely $\phi_l=l+\binom l2(t-1)+\text{higher terms}$. We deduce
\begin{equation*}\phi_l\bullet g^*\we f^*=lg^*\we f^*+\binom l2(t-1)\bullet g^*\we f^*=lg^*\we f^*-\binom l2e^*\we g^*\end{equation*}
Further $(\Lambda_k^{-1})^{\Tr}=\text{diag}(1,1,k^{-1})$ and $(\Lambda_k^{-1})^{\Tr}\we(\Lambda_k^{-1})^{\Tr}=\text{diag}(1,k^{-1},k^{-1})$. These matrices describe action of $\lambda_k$. Applying $\lambda_k$ to $\phi_l\bullet\tau$ as $\tau$ runs over the standard basis of $\X$ we complete the proof of the Lemma.\qed

The next Proposition completes the proof of the Theorem.
\begin{Prop}\label{GvsA} The set of $\g(\tl)$-orbits coincides with the set of $\A$-orbits.\end{Prop}
\pf By Corollary \ref{sizeoforbit} for every $\tau\in\X$, $\tau\g(\tl)$ is union of orbits $\tau.\omega_k\A$ for $1\le k\le p-1$. Thus it suffices to show  $\tau.\omega_k\in\tau\A$ for every $k$. We note that Lemma \ref{gaction} implies that $\omega_k\notin\A$. Nevertheless, the inclusion $\tau.\omega_k\in\tau\A$ holds for  generators $\tau$ of every orbit described in Lemmas \ref{splitorbits},\ \ref{extraorbits},\ \ref{Aorbits}. We give a sample calculation for $\tau=w_m=me^*+g^*+g^*\we f^*$ of Lemma \ref{Aorbits}. By Lemma \ref{gaction}
\begin{equation*} w_m.\omega_k=lme^*+lg^*-\binom l2e^*\we g^*+l^2g^*\we f^*\end{equation*} 
Now take $\phi$ with coordinates $u=l,v=l,r=-l^{-2}\binom l2,q=s=0$. By \eqref{Finverse} and \eqref{FwedgeF} one sees immediately that $w_m.\phi=w_m.\omega_k$.\qed
   
We can quickly dispose of the case $G=C_p\times C_p$ as promised above. Let $\tl$ denote the right regular action of $C_p$ on $R_2$. 
\begin{Prop}\rm{(\cite{Mas1})}\label{C_ptimesC_p} There are up to isomorphism $p+7$ Hopf algebras in $\mathrm{Ext}(\k^{C_p\times C_p},\k C_p)$.\end{Prop}
\pf Since all nontrivial actions of $C_p$ form a single isomorphism class we have\\ $\mathrm{Ext}(\k^{C_p\times C_p},\k C_p)=\mathrm{Ext}_{[\tl]}(\k^{C_p\times C_p},\k C_p)\cup \mathrm{Ext}_{[\t]}(\k^{C_p\times C_p},\k C_p)$. By Proposition \ref{commisotypes} $\mathrm{Ext}_{[\t]}(\k^{C_p\times C_p},\k C_p)$ contributes four nonisomorphic algebras. It remains to show that $\mathrm{Ext}_{[\tl]}(\k^{C_p\times C_p},\k C_p)$ contains $p+3$ isoclasses.

Setting $g=0$ in the definition of $G$ reduces it to $G=\Z_p\times\Z_p$ with $G\simeq R_2$ as $C_p$-module. Further reductions are as follows. The classifying space is $\X=\gen{e^*}\oplus\gen{e^*\we f^*}$,\\ $\A=\left\{\F\in\text{GL}(G)|\F=\begin{pmatrix}a_{11}&a_{12}\\0&a_{11}\end{pmatrix},\,a_{11}\ne 0\right\}$, and automorphisms $\lambda_k\in I(\tl,\tl^k)$ are defined by $\Lambda_k=\begin{pmatrix}1&0\\0&k\end{pmatrix},1\le k\le p-1$. It becomes apparent that elements $\phi$ and $\omega_k$ act on $\X$ by
\begin{align*} (a,b).\phi&=(ua,u^2b)\\(a,b).\omega_k&=(la,l^2b)\end{align*}
where $ae^*+be^*\we f^*$ is identified with $(a,b)$ and $l=k^{-1}\pmod p$ as above. Thus the orbits of $\g(\tl)$ in $\X$ coincide with those of $\A$. For the latter we note that the $\A$-orbit of every vector $(a,b)$ with $a,b\ne 0$ has $p-1$ elements, hence there are $p-1$ orbits of this kind. The set $\{(0,b)|b\ne 0\}$ is the union of two orbits, viz. $\{(0,m)|m\in{\bu{\Z_p}}^2\}$ and  $\{(0,\zeta m)|m\in{\bu{\Z_p}}^2\}$, and two more orbits $\{(0,0)\},\{(a,0)|a\ne 0\}$ are supplied by the set $\{(a,0)|a\in\Z_p\}$. \qed

We return to algebras of dimension $p^4$. We consider the case

(II) $G\simeq R_3$. We denote by $\tl_r$ the right multiplication in $R_3$. This case is sensitive to prime $p$. Let us agree to write $\mathbb{X}(\tl_r)$ as $\mathbb{X}_p$ if $G$ is a $p$-group. For $r\in\Z_pC_p$ we denote by $\ov{r}$ the image of $r$ in $R_3$. The elements $e=\ov{1},f=\ov{(t-1)},g=\ov{(t-1)^2}$ form a basis for $R_3$ in which action of $t$ is defined by $T=\begin{pmatrix}1&1&0\\0&1&1\\0&0&1\end{pmatrix}$. Let $\{e^*,f^*,g^*\}$ be the dual basis for $\G$, and $\{e^*\we f^*,e^*\we g^*,f^*\we g^*\}$ the induced basis for $\G\we\G$. We call all these bases standard. We aim to prove
\begin{Thm}\label{isoclassesAII}$\mathrm{Ext}_{[\tl_r]}(\k^{C_p^3},\k C_p)$ contains $p+9$ isoclasses, if $p>3$, and four isoclasses if $p=3$.\end{Thm}
Proof will be carried out in steps.

(1) Space $\x_p(\tl_r)$. 
\begin{Lem}\label{structureX2} If $p=3$, then
\begin{equation*} \x_3=\gen{e^*\we f^*,e^*\we g^*}\end{equation*}
For every $p>3$
\begin{equation*} \x_p=\Z_pe^*\oplus\G\we\G\end{equation*}
\end{Lem}
\pf The matrices of $t$ in the standard bases of $\G$ and $\G\we\G$ are $T^{\Tr}$ and $T^{\Tr}\we T^{\Tr}$, respectively, with $T^{\Tr}\we T^{\Tr}=\begin{pmatrix}1&0&0\\1&1&0\\1&1&1\end{pmatrix}$. From this one computes directly $(t-1)^3\bullet\G=(t-1)^3\bullet\G\we\G=0$. Since $\phi_p(t)=(t-1)^{p-1}$, it follows that $N(G)=0$ and $(\G\we\G)_N=\G\we\G$ for any $p>3$. Furhermore $\G^{C_p}=\Z_pe^*$ for every $p$. Thus as $\x_p=\G^{C_p}/N(\G)\oplus(\G\we\G)_N$ the second statement of the Lemma follows.

Say $p=3$. Then $N(\G)=(t-1)^2\bullet\G=\Z_pe^*$, hence $\G^{C_p}/N(\G)=0$. Another verification gives $(\G\we\G)_N=\gen{e^*\we f^*,e^*\we g^*}$.\qed

(2) Group $\At(\tl_r)$. For any ring $R$ with a unity viewed as a right regular $R$-module and any right $R$-module $M$ the mapping $\lambda_M:M\to\text{Hom}_R(R,M)$ defined by $\lambda_M(m)(x)=mx,x\in R$ is an $R$-isomorphism. Set $M=R=R_3$, and pick $r=a_1e+a_2f+a_3g$. The matrix of $\lambda_{R_3}(r)$ in the standard basis is $\F=\begin{pmatrix}a_1&a_2&a_3\\0&a_1&a_2\\0&0&a_1\end{pmatrix}$. From this it is evident that
\begin{equation*}\At(\tl_r)=\left\{\F\in\text{GL}(G)|\F=\begin{pmatrix}a_1&a_2&a_3\\0&a_1&a_2\\0&0&a_1\end{pmatrix}|a_i\in\Z_p,a_1\ne 0\right\}\end{equation*} 
Take $\phi=\lambda_{R_3}(r)$. Action of $\phi$ in $\G$ and $\G\we\G$ is described by $(\F^{-1})^{\Tr}$ and $(\F^{-1})^{\Tr}\we(\F^{-1})^{\Tr}$. Set $u=a_1^{-1},q=ua_2,r=ua_3$. A routine calculation gives
\begin{equation}\label{F'inverse}(\F^{-1})^{\Tr}=\begin{pmatrix}u&0&0\\-uq&u&0\\u(q^2-r)&-uq&u\end{pmatrix}\end{equation}
\begin{equation}\label{F'wedgeF'}(\F^{-1})^{\Tr}\we(\F^{-1})^{\Tr}=\begin{pmatrix}u^2&0&0\\-u^2q&u^2&0\\u^2r&-u^2q&u^2\end{pmatrix}\end{equation}
At this point it is convenient to determine a family of isomorphisms $\lambda_k:(G,\tl_r)\to(G,\tl_r^k)$. To this end, let us take $M=(R_3,\tl_r^k)$ with $2\le k\le p-1$. We set $\lambda_k=\lambda_M(e)$. By definition of $\lambda_k$ we have
\begin{equation*} \lambda_k(e)=e,\,\lambda_k(f)=e(t^k-1),\,\lambda_k(g)=e(t^k-1)^2\end{equation*}
Using the expansion $t^k-1=k(t-1)+\binom k2(t-1)^2\pmod {(t-1)^3}$ we conclude that $\Lambda_k=\begin{pmatrix}1&0&0\\0&k&\binom k2\\0&0&k^2\end{pmatrix}$ is the matrix of $\lambda_k$ in the standard basis. We shall need an explicit form of the associated matrices describing the action of $\lambda_k$ in $\G$ and $\G\we\G$, respectively. Put $l=k^{-1} \pmod p$ as usual. Then an easy calculation gives
\begin{equation}\label{asso1} (\Lambda_k^{-1})^{\Tr}=\begin{pmatrix}1&0&0\\0&l&0\\0&\binom l2&l^2\end{pmatrix},\end{equation}
\begin{equation}\label{asso2} (\Lambda_k^{-1})^{\Tr}\we(\Lambda_k^{-1})^{\Tr}=\begin{pmatrix}l&0&0\\\binom l2&l^2&0\\0&0&l^3\end{pmatrix}.\end{equation}

Unless stated otherwise we assume below that $p>3$. The degenerate case $p=3$ follows easily from the general one.

(3) Orbits of $\At(\tl_r)$ in $\x_p$. We identify an element of $\x_p$ with its coordinate vector $(a;b_1,b_2,b_3))$ relative to the standard basis of $\x_p$. We start by fixing a family of orbits separately in $\G^{C_p}$ and $\G\we\G$. These are 
\begin{align*} \Omega'_0&=\{(0)\},\Omega'_1=\{(a)|a\ne 0\},\Omega''_0=\{(0,0,0)\},\\\Omega''_{ij}&=\{(*,\ldots,*,\zeta^jb_i,0,\ldots,0)| b_i\in{\bu{\Z_p}}^2\},i=1,2,3;j=0,1\end{align*}
where the $*$ denotes an arbitrary element of $\Z_p$. For more complex orbits we need vectors $v_i(m)=(m;0,\ldots,m,0\ldots,0)$ with the second $m$ filling the $i$th slot, and $m$ running over $\bu{\Z_p}$.
\begin{Lem}\label{A'orbits} There are $3p+5$ orbits of $\At(\tl_r)$ in $\x_p$. These are
\begin{equation*}\Omega'_0\times\Omega''_0,\,\Omega'_1\times\Omega''_0,\,\Omega'_0\times\Omega''_{ij},\text{and}\;v_i(m)\At(\tl_r),i=1,2,3,j=0,1;m\in\bu{Z_p}\end{equation*}\end{Lem}
\pf It is obvious that $\Omega'_1=(1)\At(\tl_r)$. By \eqref{F'wedgeF'} $(0,\ldots,\zeta^j,0,\ldots,0).\phi\\=(*,\ldots,*,\zeta^ju^2,0\ldots,0)$ for all $i,j$ with the $*$ denoting an arbitrary element of $\Z_p$. This shows $\Omega''_{ij}=(0,\ldots,\zeta^j,0,\ldots,0)\At(\tl_r)$, hence an orbit. Similarly, by \eqref{F'inverse} and \eqref{F'wedgeF'} we have
\begin{equation}\label{v(m)phi} v_i(m).\phi=(um;*,\ldots,*,u^2m,0,\ldots,0).\end{equation}
From this one can see easily that $v_i(m)\At(\tl_r)$ has $(p-1)p^{i-1}$ elements. Another verification gives $v_i(m)\At(\tl_r)\cap v_i(n)\At(\tl_r)=\emptyset$ for $m\ne n$. Set $\Omega''_i=\Omega''_{i0}\cup\Omega''_{i1}$ and observe that $|\Omega''_i|=(p-1)p^{i-1}$ which gives $|\Omega'_1\times\Omega''_i|=(p-1)^2p^{i-1}$. Evidently $v_i(m)\in\Omega'_1\times\Omega''_i$ for all $m$ and therefore comparing cardinalities we arrive at the equality $\Omega'_1\times\Omega''_i=\bigcup_mv_i(m)\At(\tl_r)$. But clearly $\x_p=\bigcup\Omega'_k\times\Omega''_i,k=0,1;0\le i\le 3$ which completes the proof.\qed

(4) Orbits of $\g(\tl_r)$. These are listed in
\begin{Prop}\label{g'orbits} In the foregoing notation the orbits of $\g(\tl)$ in $\x_p$ are as follows.
\begin{align*} &\Omega'_0\times\Omega''_0,\,\Omega'_1\times\Omega''_0,\,\Omega'_0\times\Omega''_{1j},\,\Omega'_0\times\Omega''_2,\\
&\Omega'_0\times\Omega_{3j},\,v_1(m)\At(\tl_r),\,\Omega'_1\times\Omega''_2,\,\Omega'_1\times\Omega''_{3j}\end{align*}
where $j=0,1$ and $m$ runs over ${\Z_p}^{\bullet}$.\end{Prop}
\pf In view of  Corollary \ref{sizeoforbit} we need to determine the $\At(\tl_r)$-orbit containing $v\omega_k$ where $v$ runs over a set of generators of $\At(\tl_r)$-orbits of Lemma \ref{A'orbits}, and $\omega_k=\lambda_k\alpha_k^{-1},2\le k\le p-1$, as usual.

(i) For $\At(\tl_r)$-orbits $\Omega'_1\times\Omega''_0$ and $\Omega'_0\times\Omega''_{ij}$ generators can be chosen as \noindent$v_1=(1;0,0,0)$ and $v_{ij}=(0;0,\ldots,\zeta^j,\ldots,0)$, respectively. In view of $e^*$ and $e^*\we f^*$ being fixed points for the action of $t$, and by \eqref{asso1},  \eqref{asso2} it is immediate that 
\begin{equation}\label{v1andv1j}v_1\omega_k=lv_1\,\text{and}\, v_{1j}\omega_k=l^2v_{1j},\end{equation}
hence those sets are $\g(\tl_r)$-orbits. 

Next we take $v_{20}=e^*\we g^*$. Noting that $(t-1)^2\bullet e^*\we g^*=0$, we use the expansion $\phi_l=l+\binom l2(t-1) \pmod{(t-1)^2}$ to derive
\begin{equation*}\phi_l\bullet e^*\we g^*=c e^*\we f^*+le^*\we g^*,\,c\in\Z_p.\end{equation*}
Applying $\lambda_k$ to the last equation we find with the help from \eqref{asso2}
\begin{equation}\label{v20}e^*\we g^*\omega_k=c'e^*\we f^*+l^3e^*\we g^*,\,\text{for some}\,c'\in\Z_p.\end{equation}
The last equation shows that $v_{20}.\omega_k\in v_{21}\At(\tl_r)$ if $l$, hence $k$, is not a square, and $v_{20}.\omega_k\in v_{20}\At(\tl_r)$, otherwise. This means $v_{20}\g(\tl_r)=\Omega'_0\times(\Omega''_{20}\cup\Omega_{21}'')$, that is $\Omega'_0\times\Omega''_2$ as needed.

(ii) For a generator $v_{3j}=(0;0,0,\zeta^j)=\zeta^j f^*\we g^*$ of $\Omega_0'\times\Omega''_{3j}$ we claim that $v_{3j}\omega_k\in v_{3j}\At(\tl_r)$ for all $k$. Using the expansion $\phi_l=l+c_1(t-1)+c_2(t-1)^2 \pmod{(t-1)^3}$ we derive $\phi_l\bullet f^*\we g^*=(c_1+c_2)e^*\we f^*+c_1 e^*\we g^*+lf^*\we g^*$. Applying $\lambda_k$ we have by \eqref{asso2}
\begin{equation}\label{v30}f^*\we g^*.\omega_k=c'_1e^*\we f^*+c_1l^2e^*\we g^*+l^4f^*\we g^*,\,c'_1,c_1\in\Z_p.\end{equation}
which shows $f^*\we g^*.\omega_k\in\Omega_{3j}$. As $\Omega_{3j}=v_{3j}\At(\tl_r)$ by part (3) the claim follows.

(iii) We pause to mention that the above arguments settle the $p=3$-case. For, since $\x_3=\gen{e^*\we f^*,e^*\we g^*}$, by part (i) it has three nonzero orbits, namely $\Omega''_{1j},\Omega''_1,j=0,1$.

(iv) Here we take $v_1(m)=(m;m,0,0)$. Calculations in part (i) give $v_1(m).\omega_k=(lm;l^2m,0,0)\in v_1(m)\At(\tl_r)$ by \eqref{v(m)phi}. That is, $v_1(m)\At(\tl_r)$ is a $\g(\tl_r)$-orbit for every $m\in\bu{\Z_p}$.

It remains to show that the last three sets of the Proposition are $\g(\tl_r)$-orbits.

(v) $\Omega'_1\times\Omega''_2$ is an orbit. By Lemma \ref{A'orbits} $\Omega'_1\times\Omega''_2=\bigcup_mv_2(m)\At(\tl_r)$ where $v_2(m)=me^*+me^*\we g^*$. Note that by \eqref{v1andv1j} and \eqref{v20} there holds $v_2(m).\omega_k=(lm;c',l^3m,0)$. On the other hand we have by \eqref{v(m)phi} $v_2(n).\phi=(un;a,u^2n,0)$ where $u,a$ run over $\Z_p^{\bullet}$ and $\Z_p$, respectively. For every $l$ choosing $n=l^{-1}m, u=l^2$ and $a=c'$ we obtain $v_2(m).\omega_k\in v_2(n)\At(\tl_r)$. Letting $l$ run over $\Z_p^{\bullet}$ we see that $\bigcup_nv_2(n)\At(\tl_r)= v_2(m)\g(\tl_r)$ which completes the proof.

(vi) Here we show that each $\Omega'_1\times\Omega''_{3j}$ is an orbit. By \eqref{v20} and \eqref{v30}
\begin{equation*}v_3(m).\omega_k=(m;0,0,m).\omega_k=(lm;c',c'',ml^4)\,\text{for some}\,c',c''\in\Z_p.\end{equation*}
We seek  an $n$ such that
\begin{equation}\label{gvsA} v_3(m).\omega_k=v_3(n).\phi\,\text{for some}\,\phi\in\At(\tl_r).\end{equation}
By \eqref{v(m)phi} $v_3(n).\phi=(un;a,b,u^2n)$ with $a,b$ and $u$ taking arbitrary values in $\Z_p$ and $\Z_p^{\bullet}$, respectively. Setting $u=l^3,n=l^{-2}m, a=c'$ and $b=c''$ fullfils \eqref{gvsA}. We see $v_3(m)\g(\tl_r)=\displaystyle{\bigcup_{n\in m{\Z_p^{\bullet}}^2}}(n;0,0,n)\At(\tl_r)$. On the other hand $\Omega'_1\times\Omega''_{3j}=\displaystyle{\bigcup_{n\in m{\Z_p^{\bullet}}^2}}(n;0,0,n)\At(\tl_r)$ by comparing cardinalities of both sides.\qed

Case (B) $G=\Z_{p^2}\oplus\Z_p$ is more involved. There are $6$ classes $[\tl]$ of actions each one with its own extension theory. The final result is that $\mathrm{Ext}(\k C_p,\k^G)$ contains $3p+19$ nonisomorphic algebras $2p+7$ of which are neither commutative, nor cocommutative. The details of the proof will appear elsewhere. 
\vspace{.2in}

{\bf Appendix : Crossed product splitting of $H$}
\vspace{.1in}

{\bf Proposition}. {\em Let $H$ be an extension of $\k F$ by $\k^G$. Then $H$ is a crossed product of $\k F$ over $\k^G$}.

\pf First observe that $H$ is a Hopf-Galois extension of $\k^G$ by $\k F$ via $\rho_{\pi}=(\text{id}\otimes\pi)\Delta_H:H\to H\otimes\k F$, see e.g. the proof of \cite[3.4.3]{Mo}, hence by \cite[8.1.7]{Mo} $H$ is a strongly $F$-graded algebra. Setting $H_x=\{h\in H|\rho_{\pi}(h)=h\otimes x\}$ we have $H=\displaystyle{\oplus_{x\in F}}H_x$ with $H_1=\k^G$ and $H_xH_{x^{-1}}=\k^G$ for all $x\in F$. Next for every $a\in G$ we construct elements $u(a)\in H_x,\,v(a)\in H_{x^{-1}}$ such that 
\begin{align*} &u(a)v(a)=p_a,\,p_au(a)=u(a),v(a)p_a=v(a),\,\text{and}\\&u(a)v(b)=0\;\text{for all}\,a\ne b.\end{align*}
Indeed, were all $uv,u\in H_x,v\in H_{x^{-1}}$ lie in $\text{span}\{p_b|b\ne a\}$, then so would $H_xH_{x^{-1}}$, a contradiction. Therefore for every $a\in G$ there are $u\in H_x,v\in H_{x^{-1}}$ such that $uv=\sum c_bp_b,c_a\ne 0$. Setting $u(a)=\displaystyle\frac{1}{c_a}p_au,v(a)=vp_a$ we get elements satisfying the first three properties stated above. Furthermore, the last property also holds because $u(a)v(b)=p_au(a)v(b)p_b=p_ap_bu(a)v(b)=0$. It follows that the elements $u_x=\sum_{a\in G}u(a),v_x=\sum_{a\in G}v(a)$ satisfy $u_xv_x=1$ hence, as $H$ is finite-dimensional, $v_xu_x=1$ as well. Thus $u_x$ is a $2$-sided unit in $H_x$. 

Now define $\chi:\k F\to H$ by $\chi(x)=\displaystyle\frac{1}{\epsilon_H(u_x)}u_x$. One can see immediately that $\chi$ is a convolution invertible mapping satisfying $\rho_{\pi}\circ\chi=\chi\otimes\text{id},\chi(1_F)=1$ and $\epsilon_H\circ\chi=\epsilon_F$. Thus $\chi$ is a section of $\k F$ in $H$, which completes the proof.\qed


\begin{thebibliography}{999}

\bibitem{A} N. Andruskiewitsch, Notes on extensions of Hopf algebras, {\em Can. J. Math.} {\bf 48}(1)(1996), 3-42.
\bibitem{BT} F. R. Beyl and J. Tappe, Group Extensions, Representations, and the Schur Multiplicator, Lecture Notes in Mathematics {\bf 958}, Springer-Verlag, 1982.
\bibitem{BCM} R.J. Blatttner, M. Cohen and S. Montgomery, Crossed Product and Inner Actions of Hopf Algebras, {\em Trans. Amer. Math. Soc} {\bf 298}(2)(1986),671-711.
\bibitem{Bou} N. Bourbaki, Elements of Mathematics, Algebra I, Springer-Verlag, 1989.
\bibitem{B} N.P. Byott, Cleft extensions of Hopf algebras, {\em J. Algebra} {\bf 157}(1993),405-429.
\bibitem{Ha} M. Hall, Jr., ``The Theory of Groups'', The Macmillan Company, New York, 1959. 
\bibitem{H} I. Hofstetter, Extensions of Hopf algebras and their cohomological description, {\em J. Algebra} {\bf 164}(1994), 264-298.
\bibitem{KMM} Y. Kashina, G. Mason, S. Montgomery, Computing the Frobenius-Schur indicator for abelian extensions of Hopf algebras, {\em J. Algebra} {\bf 251}(2002), 888-913.
\bibitem{La} S. Lang, ``Algebra'', Addison-Wesley, 1993.
\bibitem{LR} R.G. Larson and D.E. Radford, Semisimple cosemisimple Hopf algebras, {\em Amer. J. Math.} {\bf 110}(1988), 187-195.
\bibitem{Ma} M. Mastnak, Hopf algebra extensions arising from semi-direct products of groups, {\em J. Algebra} {\bf 251}(2002), 413-434.
\bibitem{Mac} S. MacLane, ``Homology'', Die Grundlehren der Mathematischen Wissenschaften {\bf 114}, Springer-Verlag, 1963.
\bibitem{MD} A. Masuoka and Y. Doi, Generalization of cleft comodule algebras, {\em Comm. Algebra} {\bf 20}(1992), 3703-3721.
\bibitem{Mas1} A. Masuoka, Self-dual Hopf algebras of dimension $p^3$ obtained by extensions, {\em J. Algebra} {\bf 178}(1995), 791-806.
\bibitem{Mas2} A.Masuoka, The $p^n$ theorem for semisimple Hopf algebras, {\em Proc. Amer. Math. Soc.},{\bf 124(3)}(1996), 735-737.
\bibitem{M} A. Masuoka, Extensions of Hopf algebras (Lecture Notes, University of Cordoba, 1997) Notas Mat.No. 41/99, FaMAF Uni. Nacional de Cordoba, 1999.
\bibitem{Mo} S. Montgomery, Hopf Algebras and their Actions on Rings, in: {\em CMBS Reg. Conf. Ser.Math.} {\bf 82}, AMS, 1993.
\bibitem{NVO} C. Nastasescu and F. Van Oystaeyen, On strongly graded rings and crossed products, {\em Comm. Algebra} {\bf 10} (1982), 2085-2106.
\bibitem{S1} H.-J. Schneider, Some remarks on exact sequences of quantum groups, {\em Comm. Algebra} (9){\bf 21}(1993), 3337-3358.
\bibitem{S2} H.-J. Schneider, A normal basis and transitivity of crossed products for Hopf algebras, {\em J. Algebra} {\bf 152}(1992), 289-312.
\bibitem{St} D. Stefan, The set of Types of $n$-dimensional semisimple and cosemisimple Hopf algebras is finite, {\em J. Algebra} {\bf 193}(1997), 571-580.
\bibitem{T} M. Takeuchi, Matched pairs of groups and bismash products of Hopf algebras, {\em Comm. Algebra} {\bf 9}(1981), 841-882.
\bibitem{Y} K. Yamazaki, On projective representations and ring extensions of finite groups, {\em J. Fac. Science Univ. Tokyo, Sect. I} {\bf 10}(1964), 147-195.














\end{thebibliography}
\end{document}